\documentclass[12pt]{article}
\usepackage{amsmath,amssymb}
\usepackage[mathscr]{euscript}
\usepackage{color}
\usepackage{epsfig}
\usepackage[all]{xypic}

\newtheorem{prop}{Proposition}[section]
\newtheorem{rema}{Remark}[section]
\newtheorem{defi}{Definition}[section]
\newtheorem{lemm}{Lemma}[section]
\newtheorem{theo}{Theorem}[section]
\newtheorem{coro}{Corollary}[section]

\newcommand{\bbox}{\normalsize {}%
        \nolinebreak \hfill $\blacksquare$ \medbreak \par}

\newcommand{\R}{\mathbb{R}}
\newcommand{\C}{\mathbb{C}}

\newcommand{\N}{\mathbb{N}}

\newcommand{\T}{\mathbb{T}}
\newcommand{\F}{\mathbb{F}}

\newcommand{\caus}{\hbox{\tiny{Cau}}^s}

\newcommand{\Caus}{\hbox{Cau}^s}

\newcommand{\lc}{\boldsymbol{[}\!\!\boldsymbol{[}}
\newcommand{\rf}{\boldsymbol{]}\!\!\boldsymbol{]}}

\newcommand{\bphi}{\boldsymbol{\phi}}

\newcommand{\dirac}{\partial\!\!\!/}

\usepackage{stmaryrd}                       

\title{First integrals for nonlinear dispersive equations}
\author{Fr{\'e}d{\'e}ric \textsc{H{\'e}lein}\footnote{Institut de Math{\'e}matiques de Jussieu---Paris Rive Gauche,
UMR CNRS 7586, Universit{\'e} Paris Diderot --- Paris 7,
B{\^a}timent Sophie Germain, Case 7012, 
75205 Paris Cedex 13, France, \textsf{helein@math.univ-paris-diderot.fr}}
}

\begin{document}
\maketitle
\begin{abstract}
\emph{Given a solution of a semilinear dispersive partial differential equation with a real analytic nonlinearity,
we relate its Cauchy data at two different times by nonlinear representation formulas in terms of convergent
series. These series are constructed by means of generating functions. All this theory is based on a
new suitable formulation of the dynamics of solutions of dispersive equations.}
\end{abstract}

Consider a partial differential equation
\begin{equation}\label{equa0}
 Lu + N(u,\partial u) = 0,
\end{equation}
which describes the evolution of a map $u$ from a space-time $\R^{n+1}$ (with coordinates
$(x^0,\cdots ,x^n)$) to a finite dimensional vector space.
Here $L$ is a linear differential operator (e.g. the wave operator $\square
= \partial_0^2 - \Delta$, Klein--Gordon $\square + m^2$, Dirac
$\dirac + im$, or any combination) and $N$ is a real analytic nonlinear function on $u$
and its first space-time derivatives $\partial u$.
For any $t\in \R$, denote by $[u]_t$ the Cauchy data of $u$ at time $t$. We
address the question: assume that we know $[u]_{t_1}$ for some $t_1\in \R$,
can we compute the value of $u$ at a point at another time $t_2$ ?
If $N$ is a linear function the answer is positive and is given by a linear integral formula,
if $N$ is a polynomial this may also work by using series, i.e. an infinite sum of
multilinear integrals, as we will present here. In \cite{harrivel} D. Harrivel
obtained such a result for a (roughly speaking $\mathcal{C}^2$)
solution of the scalar Klein--Gordon equation $\square u + m^2u + \lambda u^2 = 0$.
It amounts to build a time dependant family of functionals $(\mathcal{S}_t)_t$
of Cauchy data s.t., if $u$ solves (\ref{equa0}), then $\mathcal{S}_t([u]_t)$ does not depend on $t$. Moreover
one can prescribe $\mathcal{S}_{t_2}$ to be any linear functional. 
By choosing e.g. $\mathcal{S}_{t_2}$ to be the Dirac distribution at some point
we thus get a positive answer of the previous question. The functionals 
$\mathcal{S}_t$ are series, each term of
which is a sum of integrals over Cartesian products of the space-time
built from planar binary trees by using Feynman rules. The important
point is that one can ensures that the series converges for $|t_2-t_1|$
sufficiently small.

In \cite{harrivel-helein} this result
was further extended to systems with more general (real analytic) nonlinearities and
for less regular solutions (roughly speaking $\mathcal{C}^1$). The method, which  was
different from \cite{harrivel}, did not use a combinatorial
analysis of the series, but rests on the construction of a generating
function which, by using Wick's theorem for developping it, gives us the
desired expansion.

The following paper presents an improvement of the results in \cite{harrivel-helein}.
A new ingredient is a different formulation of the dynamics,
which allows us to deal with even less regular solutions (roughly speaking $\mathcal{C}^0$
in general).
This formulation is, we believe, new although it is a straightforward consequence
of the well-known Duhamel formula. To explain it, consider the
standard way to formulate an evolution PDE such as (\ref{equa0}) as an
ODE in the infinite dimensional space of all Cauchy data:
\begin{equation}\label{equa-wrong-way}
 \frac{d[u]_t}{dt} = X([u]_t).
\end{equation}
We introduce an alternative formulation: we work in the space $\mathcal{E}_0$ of solutions to 
Equation (\ref{equa-linear}) below and replace $[u]_t$ by the unique
solution $\varphi$ to the linear equation
\begin{equation}\label{equa-linear}
 L\varphi = 0
\end{equation}
the \emph{Cauchy data of which is $[u]_t$}, i.e. the same as $u$ at time $t$. We denote by
$\Theta_tu\in \mathcal{E}_0$ this solution. Let $G$ be the homogeneous solution to $LG=0$
s.t., if $Y:M\longrightarrow \R$ is the function defined by $Y(x) = 1$ if $x^0\geq 0$ and $G(x) = 0$ if
$x^0<0$, for some time coordinate $x^0$, then $YG$ is the retarded fundamental solution of $L$.
Consider the time dependant vector field $(V_t)_t$ on $\mathcal{E}_0$ defined by
$V_t\varphi:= \int_{y^0=t}d\vec{y}G(\cdot -y)N(\varphi,\partial \varphi)(y)$.
Our first result is:
\begin{theo}\label{theo-zero}
A map $u$ is a solution of (\ref{equa0}) if and only if the map 
$t\longmapsto \Theta_tu$ is a solution to
\begin{equation}\label{equa-right-way}
 \frac{d\Theta_tu}{dt} + V_t(\Theta_tu) = 0.
\end{equation}
\end{theo}
A precise statement of this result is the content of Theorem \ref{theo-magic-dynamics}.
An advantage of Equation (\ref{equa-right-way}) is that it is manifestly covariant: the
space $\mathcal{E}_0$ in which $\Theta_tu$ takes values does not depend on $t$ nor on
any choice of space-time coordinates, in contrast with the target space of $t\longmapsto [u]_t$.
This advantage is even more striking on a curved space-time, where a similar result
will be proved (Theorem \ref{theo-magic-dynamics-curved}).
A second advantage is that the map $t\longmapsto \Theta_tu$
is more regular than $t\longmapsto [u]_t$: under general hypotheses, if $u$ is a weak
solution of (\ref{equa0}) then $t\longmapsto \Theta_tu$ is $\mathcal{C}^1$~!

This formulation is also useful for the problem expounded previously.
Consider the space $\mathbb{F}$ of real analytic functionals on $\mathcal{E}_0$.
We define for all $t$ the first order `differential' linear operator $V_t\cdot$
acting on $\mathbb{F}$ by:
\begin{equation}\label{define-bdelta-0}
\forall f\in \F,
 \forall \varphi\in \mathcal{E}_0,\quad (V_t\cdot f)(\varphi):= \delta f_\varphi(V_t(\varphi)),
\end{equation}
where $\delta f_\varphi$ is the differential of $f$ at $\varphi$.
Then one of our main result is that we can make sense of the chronological exponential
$T\hbox{exp}\left(\int_{t_1}^{t_2} dsV_s\cdot \right)$ as
a linear operator acting on $\mathbb{F}$, continuous in a suitable topology, if
$|t_2-t_1|$ is sufficiently small.
This operator is the key for constructing the family $(\mathcal{S}_t)_t$ of
operators such that $\mathcal{S}_t([u]_t)$ does not depend on $t$ if $u$ is a solution
of (\ref{equa0}): 
\begin{theo}\label{theo-intro}
Let $r>0$.
 There exists a constant $\overline{t}>0$ which depends on Equation (\ref{equa0}) and on
 $r$, such that, for any $t_1,t_2\in \R$ such that $|t_2-t_1|<\overline{t}$ and for any
 $f\in \mathbb{F}$, with a radius of convergence $r$, the functional
 \[
  U_{t_1}^{t_2}f:= T\hbox{exp}\left(\int_{t_1}^{t_2} dsV_s\cdot \right)f
 \]
is well defined on a ball in $\mathbb{F}$ and has a non vanishing radius of convergence $R$.
Moreover, if $u$ is a solution of (\ref{equa0}) the Cauchy data of which is smaller than $R$, then
\begin{equation}\label{main-quantity}
(U_{t_1}^{t_2}f)(\Theta_{t_2}u) \quad
\hbox{is equal to }f\left(\Theta_{t_1}u\right).
\end{equation}
\end{theo}
Details on the statement in Theorem \ref{theo-intro} (the topology on
$\mathcal{E}_0$ and on the space of Cauchy data) will given in the next Section. In general
we will set $u\in \mathcal{C}^0(I,H^s(\mathbb{R}^n))\cap \mathcal{C}^1(I,H^{s-r}(\mathbb{R}^n))$,
where $r$ depend on $L$ (e.g. $r=1$ for $L= \square$) and  $s>n/2$ in general.
However for a Klein--Gordon equation with some polynomial nonlinearity, it may work
for some special values of $s$ and $n$ s.t. $s\leq n/2$ (see Remark \ref{remarksursegal1}).

This result can be restated in a different language inspired by perturbative quantum fields theory:
$U_{t_1}^{t_2}f$ can be written
\[
 U_{t_1}^{t_2}f = \left(T\hbox{exp}\int_{t_1<y^0<t_2}dy\
 N^i(\bphi,\partial \bphi)(y)\bphi^+_i(y)\right)f,
\]
where $\bphi$ and $\bphi^+$ are kind of \emph{creation} and \emph{annihilation} operators respectively
(see Section \ref{sec-comparison} for details).\\

\noindent
\textbf{Plan of the paper}

For simplicity most results are presented for a differential operator with constant coefficients
on a flat space-time.
Section 1 contains the notations and a precise formulation of the hypotheses
needed for the theory on a flat space-time.
In Section 2 we construct the map $u\longmapsto \Theta_tu$ and the vector field $V_t$
on a flat space-time.
We end with the proof of Theorem \ref{theo-magic-dynamics}, a version of Theorem \ref{theo-zero} on a flat
space-time. We also show that $V_t$ is real analytic
on an open ball in $\mathcal{E}_0$.
In Section 3 we extend these results to a curved space-time.
For simplicity we restrict ourself to the Klein--Gordon operator and a cubic
nonlinearity. We show also that the formulation (\ref{equa-right-way}) works
if we replace a foliation by space-like hypersurfaces which are the level sets of
a time function by a more general family of space-like hypersurfaces which may overlap.

In Sections 4, 5 and 6 we developp a theory valid in any Banach space $\mathbb{X}$.
In Section 4 we introduce various topologies on the space $\mathbb{F}(\mathbb{X})$ of
real analytic functions on bounded balls of $\mathbb{X}$.
We define in particular, for any $r\in [0,+\infty]$, the space $\mathbb{F}_r(\mathbb{X})$ of real
analytic functions on $\mathbb{X}$ which, roughly speaking, have a radius of convergence
greater or equal to $r$. We also
derive properties satisfied by a time dependant family $(V_t\cdot)_t$ of real analytic first order
differential operators acting on $\mathbb{F}(\mathbb{X})$.
In Section 5 we prove
the existence of the chronological exponential $U_{t_1}^{t_2} = T\hbox{exp}\left(\int_{t_1}^{t_2} dsV_s\cdot\right)$
as a bounded operator acting between subspaces of $\mathbb{F}(\mathbb{X})$, if $|t_2-t_1|$ is small enough.
The difficulty is that
the operators $V_t\cdot$ are not bounded in any topology. Hence the chronological
exponential cannot make sense as a bounded operator from a topological to itself.
However we will prove that $U_{t_1}^{t_2}$ maps
continuously $\mathbb{F}_R(\mathbb{X})$ to $\mathbb{F}_{e^{-|t_2-t_1|X}(R)}(\mathbb{X})$,
where $X$ is a (positive) real analytic vector field on $\R$ which is constructed out of
Equation (\ref{equa0}) and of the choice of topology on the set of its solutions.
In Section 6 we prove that $(U_{t_1}^{t}f)(\varphi(t))$ does not depend on $t$ if
$\varphi$ is a solution of $\frac{d\varphi}{dt}+V_t(\varphi)=0$, a result which,
combined with Theorem \ref{theo-magic-dynamics} or Theorem \ref{theo-magic-dynamics-curved},
implies different versions of Theorem \ref{theo-intro}.

In Section 7 and 8 we give some applications of our results and discuss the analogy and
the difference with methods from Quantum Field Theory.\\

\noindent
\textbf{Further comments}

This work is motivated by questions in \cite{helein,heleinLeeds}.
Our formulation of the dynamics by (\ref{magic-dynamics}) can be viewed as an analogue for dispersive partial differential
equations of Lagrange's method of \emph{variation of the constant}, it is also an infinitesimal version of Duhamel's formula
(\ref{duhamelforever}). This is the reason for the name `Lagrange--Duhamel' for $V_t$.

Developping (\ref{main-quantity}) by using Wick theorem leads to an expansion
in terms of `Feynman trees', as for instance (\ref{formuleexemple}).
A heuristic way to understand where this comes from consists in inserting the l.h.s. of
(\ref{formule2.6}) in the integral in the r.h.s. of it and in iterating this process. Then
we see easily that $u(x)$ should be expressed as the sum of a formal series. But
it seems difficult to prove directly by this method that this process converges and to 
estimate the radius of convergence of the series. On the other hand this process is also
the key of the Picard fixed point Theorem which is used to prove the local existence of solutions.
However the proof of the fixed point Theorem is based on precise estimates of the previous
process but it hides the structure of the series which is generated by this process. Our result
can hence be understood as filling the gap between both ways.

Series expansions of solution to nonlinear \emph{ordinary} differential equations (ODE) have a long history.
We can mention Lie series
defined by K.T. Chen \cite{chen}, the Chen--Fliess series \cite{fliess} introduced in the framework of control
theory by M. Fliess (or some variants like Volterra series or Magnus expansion \cite{magnus}) which are extensively
used in control theory \cite{agrachev-gamkrelidze,sussman,kawski-sussman} but also in the study of dynamical systems
and in numerical analysis. Other major tools are Butcher series which explain the structure of Runge--Kutta methods
of approximation of the solution of an ODE. They have been introduced by J.C. Butcher \cite{butcher} and developped
by E. Hairer and G. Wanner \cite{hairerwanner} which explain that Runge--Kutta methods are gouverned by trees. Later on
C. Brouder \cite{brouder,brouder2} realized that the structure which underlies the original Butcher's computation 
coincides with the Hopf algebra defined by D. Kreimer in his work about the renormalization theory \cite{kreimer}.
Concerning analogous results on nonlinear \emph{partial} differential equations,
it seems that the fact that one can represent solutions or functionals
on the set of solutions by series indexed by trees is known to physicists since the work of J. Schwinger and R. Feynman
(and Butcher was also aware of that in his original work). However it is not that easy to find precise references
in the litterature: the Reader may look e.g. at  \cite{duetsch}, where a formal series expansion is presented and the recent
paper \cite{finster} for comparison with quantum field theory.
But, to our knowledge, the only rigorous results (i.e. with a proof of convergence of the series)
can be found in \cite{harrivel,harrivel-helein}.

We have used relatively elementary tools from the analysis of PDE's and, in particular, we do not rely on the modern
theory for wave and Schr{\"o}dinger equations (Strichartz estimates, Klainerman bilinear estimates, etc.). Many
improvements in these directions could be provided, although they may not be straightforward.
Also we are not able to apply our theory
the KdV equation, since its nonlinearity cannot be controlled by our methods. 
Another question concerns the extension of our results to an infinite time
interval and to relate together the asymptotic data for $t\rightarrow -\infty$ and $t\rightarrow + \infty$.
One may indeed ask whether the limits $u_\pm:= \lim_{t\rightarrow \pm\infty}\Theta_tu$ exist and, if so, if
for $f\in \mathbb{F}$, 
\[
f(u_-) =
 \left(\left(T\hbox{exp}\int_{-\infty}^0 d\tau V_\tau\cdot\right)f\right)(\Theta_0u) =
 \left(\left(T\hbox{exp}\int_{-\infty}^{+\infty} d\tau V_\tau\cdot\right)f\right)(u_+).
\]
Such identities (and their analogues by exchanging $u_+$ and $u_-$) would imply in particular that
the scattering map $S:u_-\longmapsto u_+$ and the wave operators $W_\pm:u_\pm\longmapsto \Theta_0u$ are well-defined
and real analytic\footnote{In our definition $S$ and $W_\pm$ map $\mathcal{E}^s_0$ to itself. This differs
from most references where the scattering map reads in our notations $\Phi_0^{-1} \circ S\circ \Phi_0: \Caus
\longrightarrow \Caus$ and the wave maps are $\Phi_0^{-1} \circ W_\pm\circ \Phi_0: \Caus
\longrightarrow \Caus$ (see Paragraph \ref{paragraph2.2} for the definition of $\Phi_0$).}.
In the light of results in \cite{morawetz-strauss,strauss,strauss1,raczka-strauss,brenner,baez,baez-zhou}
this should be true for the equation $\square u + u^3 = 0$ on $\mathbb{R}^{1+3}$ and for $s=1$,
due to dispersive effects (Strichartz estimates).
The key point in all these works is an estimate of the type
$\int_{-\infty}^{+\infty}\left(\int_{\mathbb{R}^3}u^6d\vec{x}\right)^{1/2}dt<+\infty$,
which, e.g., holds for a solution $u$ of $\square u + u^3 = 0$ with finite energy.\\

\noindent
\emph{Acknowledgements} --- I wish to thank Isabelle Gallagher for explanations about paraproducts.
This paper is a extended and improved version of an earlier work in collaboration with
Dikanaina Harrivel \cite{harrivel-helein}.

\section{Notations and hypotheses}\label{notationshypotheses}
\textbf{Generalities} ---
$M:= \mathbb{R}\times \mathbb{R}^n$ represents an $(n+1)$-dimensional flat space-time.
We denote by $x = (x^0,\vec{x}) = (x^0,x^1,\cdots, x^n)$ the coordinates on $M$ and set
$\partial_\mu = \frac{\partial}{\partial x^\mu}$ for $0\leq \mu\leq n$.
We let $E$ be a finite dimensional real vector space and we consider maps from $M$ to $E$.

For any smooth fastly decreasing functions $f\in{\cal S}(\R^n)$ we
define its Fourier transform $\hat{f}(\xi) = \int_{\R^n} f(\vec{x})e^{-i\vec{x}\cdot \xi}d\vec{x}$
and we extend it to space ${\cal S}'(\R^n)$ of tempered distributions by the standard
duality argument. In case of a map $f$ which depends on $(t,\vec{x})\in I\times \mathbb{R}^n$, we also denote by
$\hat{f}(t,\xi) = \int_{\R^n} f(t,\vec{x})e^{-i\vec{x}\cdot \xi}d\vec{x}$ the partial Fourier transform
in space variables. 

For $s\in \R$, we let
$H^s(\R^n):= \{\varphi\in {\cal S}'(\R^n)|\ [\xi\longmapsto
\langle\xi\rangle^s\widehat{\varphi}(\xi)]\in L^2(\R^n)\}$,
where $\langle\xi\rangle:=\sqrt{m^2+|\xi|^2}$ and we set
$||\varphi||_{H^s}:= ||\langle\xi\rangle^s\widehat{\varphi}||_{L^2}$.
We let $H^s(\mathbb{R}^n,E)$ be the Sobolev space of $E$-valued maps on $\mathbb{R}^n$.
If $\varphi\in H^s(\R^n,E)$ has the coordinates $\varphi^i$
($1\leq i\leq \hbox{dim}E$) in a basis of $E$ we set
\begin{equation}\label{def-HsRnE}
 \|\varphi\|_{H^s}:= \sum_{i=1}^{\hbox{\footnotesize{dim}}E} \|\varphi^i\|_{H^s}.
\end{equation}

\noindent
\textbf{The class of differential operators $L$} ---
We suppose that there is a splitting $E:= E_1\oplus E_2$, where $E_1$ and $E_2$ 
are two vector subspaces of $E$.
This leads to a decomposition of any map $\varphi:M\longrightarrow E$ 
as $\varphi = \left(\varphi_1, \varphi_2\right)$.
We assume that the linear differential operator $L$ acting on smooth
maps $u:M\longrightarrow E$ has the form
\begin{equation}\label{defL}
L = \left(\begin{array}{cc}L_1 & 0 \\ 0 & L_2\end{array}\right)
= \left(\begin{array}{cc}\gamma^0\partial_0 + P_1(\vec{\partial}) & 0 \\
0 & \partial_0^2 + P_2(\vec{\partial})\end{array}\right),
\end{equation}
where $\gamma^0\in \hbox{End}(E_1)$ is an invertible matrix,
$\vec{\partial}:= (\partial_1,\cdots ,\partial_n)$ and $P_1$ and $P_2$ are polynomials
with coefficients in respectively $\hbox{End}(E_1)$ and $\hbox{End}(E_2)$ and of degree respectively
$r$ and $2r$, where $r\in \mathbb{N}^*$.

We assume that, $\forall \xi\in (\R^n)^*$, $i(\gamma^0)^{-1}P_1(i\xi)$ is a Hermitian matrix
and $P_2(i\xi)$ is positive Hermitian. Moreover we suppose that there exists
two constants $\alpha>0$ and ${\mu_0}\geq 0$ s.t., in the sense of Hermitian operators acting on $E_2$,
\begin{equation}\label{massive-nonmassive}
\forall \xi\in (\mathbb{R}^n)^*,\quad
P_2(i\xi) \geq \alpha({\mu_0} + |\xi|^{2r}).
\end{equation}
Below is a list of examples for $L$ (setting $\square = \partial_0^2 - \Delta$).
\[
 \begin{array}{|cc|c|c|c|c|c|}
\hline  \hbox{} & L & E_1 & E_2 & r\\
\hline  \hbox{Klein--Gordon} & \square + m^2 & \{0\} & \mathbb{R} &  1\\
\hline  \hbox{Schr{\"o}dinger} & i\partial_0
+\Delta & \mathbb{C} & \{0\} & 2\\
\hline  \hbox{Dirac on }\mathbb{R}^4
& \dirac + im = \gamma^\mu\partial_\mu + im & \mathbb{C}^4 & \{0\} & 1\\
\hline   \begin{array}{c}
          \hbox{linearized}\\
\hbox{Korteweg--de Vries}
         \end{array} & \partial_0
+ (\partial_1)^3 & \mathbb{R} & \{0\} & 3\\
\hline  \begin{array}{c}
         \hbox{linearized} \\ \hbox{Dirac--Maxwell}\\ \hbox{(in Lorentz gauge)}
        \end{array} &  \left(\begin{array}{cc}
           \dirac + im & 0 \\
              0 &  \square
           \end{array}\right) & \mathbb{C}^4 &  \mathbb{R}^4 &  1\\
\hline
 \end{array}
\]

\noindent
\textbf{The function spaces} --- For any $s\in \mathbb{R}$ and any interval $I\subset \R$ we define the space
\[
 \mathcal{F}^s(I):= \mathcal{C}^0(I,H^s(\mathbb{R}^n,E))
\cap \mathcal{C}^1(I,H^{s-r}(\mathbb{R}^n,E))
\]
on which the operator $L$ acts.
The natural space of Cauchy data for $L$ on $\mathcal{F}^s(I)$ is
$\hbox{Cau}^s :=  H^s(\mathbb{R}^n,E)\times H^{s-r}(\mathbb{R}^n,E_2)$.
For any $(\psi,\chi)\in \hbox{Cau}^s$, we set
\[
 \|\psi, \chi\|_{\caus}:= \|\psi\|_{H^s} + \|\chi\|_{H^{s-r}}.
\]
The space $\mathcal{F}^s(I)$ is equipped with the norm
$\|u\|_{\mathcal{F}^s(I)}:= \sup_{\tau\in I}\|[u]_\tau\|_{\caus}
= \|u\|_{L^\infty(I,H^s)} + \|\partial_0u_2\|_{L^\infty(I, H^{s-r})}$.

For any map $\varphi$ defined on a neighbourhood of $\{t\}\times \mathbb{R}^n$
in $M$, define its Cauchy data at time $t$ by
$[\varphi]_t:= (\varphi|_t,\partial_0 \varphi_2|_t) \in \hbox{Cau}^s$,
where, for any function $\psi$, we note $\psi|_t:= \psi(t,\cdot)$ its restriction 
to $\{t\}\times \mathbb{R}^n$ (which we identify with a function defined on $\mathbb{R}^n$).
For any $I\subset \R$ and $t\in I$, this hence defines a continuous linear map of norm one
\begin{equation}\label{trivialmap}
 \begin{array}{ccc}
  \mathcal{F}^s(I) & \longrightarrow & \hbox{Cau}^s\\
  \varphi & \longmapsto &[\varphi]_t
 \end{array}
\end{equation}
For any interval $I\subset \R$ we define the space of solutions to the linear equation $L\varphi=0$:
\begin{equation}\label{definition-E0s}
 \mathcal{E}_0^s(I):=\{ \varphi\in \mathcal{F}^s(I)|\
L\varphi = 0\hbox{ weakly}\}.
\end{equation}
This space is equipped with the norm $\|u\|_{\mathcal{F}^s(I)}$.

By Proposition \ref{proposition-continuity-of-Phi},
assuming Hypotheses (\ref{defL}) and (\ref{massive-nonmassive}), for any $t\in \mathbb{R}$ and any pair
$(\psi,\chi) \in \hbox{Cau}^s$, there exists an unique map $\varphi\in \mathcal{E}_0^s(\R)$
s.t. $[\varphi]_t = (\psi, \chi)$, i.e. a solution $\varphi\in \mathcal{F}^s(\R)$ of:
\begin{equation}\label{linear-pb}
L\varphi = 0\quad\hbox{s.t.}\quad \varphi|_t = \psi\quad
\hbox{and}\quad \partial_0\varphi|_t = \chi.
\end{equation}
We denote by $\Phi_t(\psi,\chi)$ this solution.\\

\noindent
\textbf{The map $\Theta$} --- For any map $u$ defined on a neighbourhood of $\{t\}\times \mathbb{R}^n$
we set
\[
\Theta_tu:= \Phi_t\left([u]_t\right), 
\]
i.e. $\Theta_tu$ is the unique solution $\varphi$ of $L\varphi=0$ s.t. $[\varphi]_t = [u]_t$. This hence defines the map
\begin{equation}\label{composed-maps}
 \begin{array}{cccccc}
  \Theta: & I\times \mathcal{F}^s(I) & \longrightarrow &
I\times \hbox{Cau}^s & \longrightarrow & \mathcal{E}_0^s(\R)
\\
 & (t,u) & \longmapsto & (t,[u]_t) & \longmapsto & \Theta_tu
 \end{array}
\end{equation}

\noindent
\textbf{Polynomials and real analytic functions} --- Let $\mathbb{X},\mathbb{Y}$ be two Banach
spaces and $p\in \mathbb{N}$. For any $r>0$, denote by $B_\mathbb{X}(r)$ the open
ball of radius $r$ and of center $0$ in $\mathbb{X}$.
A linear map $f_\otimes $ from $\mathbb{X}^{\otimes p}$ to $\mathbb{Y}$ is \emph{symmetric} if
$\forall \varphi_1,\cdots, \varphi_p\in \mathbb{X}$, $\forall \sigma\in \mathfrak{S}(p)$,
$f_\otimes (\varphi_{\sigma(1)}\otimes  \cdots \otimes  \varphi_{\sigma(p)}) =
f_\otimes (\varphi_1\otimes  \cdots\otimes  \varphi_p)$, where 
$\mathfrak{S}(p)$ is the symmetric group with $p$ elements.
A \emph{homogeneous polynomial map $f:\mathbb{X}\longrightarrow \mathbb{Y}$ of degree $p$} is a map such that there exists a
symmetric linear map $f_\otimes:\mathbb{X}^{\otimes p} \longrightarrow \mathbb{Y}$ such that $\forall \varphi \in \mathbb{X}$,
$f(\varphi) = f_\otimes (\underbrace{\varphi\otimes \cdots \otimes \varphi}_p)$. Note that
$f_\otimes$, if it exists, is unique and is given by the polarization formula:
\begin{equation}\label{polarization}
 f_\otimes (\varphi_1\otimes \cdots \otimes \varphi_p) =
\frac{1}{2^pp!}\sum_{\epsilon = (\epsilon_1,\cdots,\epsilon_p)\in \{\pm 1\}^p}
\left(\prod_{j=1}^p\epsilon_j\right) f\left(\sum_{j=1}^p\epsilon_j \varphi_j\right).
\end{equation}
If so we denote by $\|f\|_{\otimes}$ the smallest nonnegative constant
such that $\forall \varphi_1,\cdots, \varphi_p\in \mathbb{X}$,
\begin{equation}\label{def-tensornorm}
 \|f_\otimes (\varphi_1\otimes \cdots \otimes \varphi_p)\|_\mathbb{Y}
\leq \|f\|_{\otimes} \|\varphi_1\|_\mathbb{X} \cdots \|\varphi_p\|_\mathbb{X}
\end{equation}
Most of the time we will abuse notations identifying $f_\otimes $ with $f$, when there is no ambiguity.
We denote by $\mathcal{Q}^p(\mathbb{X},\mathbb{Y})$ the vector space of homogeneous polynomial maps
from $\mathbb{X}$ to $\mathbb{Y}$ of degree $p$.

A \emph{formal series $f$ from} $\mathbb{X}$ to $\mathbb{Y}$ is an infinite sum
\begin{equation}\label{F-splitting}
  f = \sum_{p=0}^\infty f^{(p)},
\end{equation}
where $\forall p\in \N$, $f^{(p)}\in \mathcal{Q}^p(\mathbb{X},\mathbb{Y})$.
The \emph{multiradius of convergence}\footnote{Note that
beside $\|f^{(p)}\|_\otimes$ defined by (\ref{def-tensornorm}), one can also
define $\|f^{(p)}\|:= \inf_{\varphi\in \mathbb{X}\setminus \{0\}}
\|f^{(p)}(\varphi)\|_\mathbb{Y}/\|X\|_\mathbb{X}^p$ and the radius of convergence
$\rho(f)$ of the series $\sum_{p=0}^\infty \|f^{(p)}\|z^p$.
One can then prove by using (\ref{polarization})
that $\|f^{(p)}\|\leq \|f^{(p)}\|_{\otimes} \leq \frac{p^p}{p!}\|f^{(p)}\|$,
which implies by using Stirling's formula
that $e^{-1}\rho(f)\leq \rho_\otimes(f) \leq \rho(f)$.} of $f$ is the radius of convergence of the series
\begin{equation}\label{formule-qui-donne-VexX}
\lc f \rf (z):= \sum_{p=0}^\infty \|f^{(p)}\|_{\otimes} z^p
\end{equation}
and is denoted by $\rho_\otimes(f)$.
We denote by $\mathbb{F}(\mathbb{X},\mathbb{Y})$ the space of formal series from
$\mathbb{X}$ to $\mathbb{Y}$. If $\rho_\otimes(f)>0$, $f$ defines
a \emph{real analytic map} from $B_\mathbb{X}(\rho_\otimes(f))$ to $\mathbb{Y}$ by the relation
$\forall \varphi\in B_\mathbb{X}(\rho_\otimes(f))$, $f(\varphi) = \sum_{p=0}^\infty f^{(p)}(\varphi)$.
This map is continuous (Lemma \ref{lemma-series-are-continuous}) and
satisfies the inequality
\begin{equation}\label{inequality-VcontroleparVexX}
 \|f(\varphi)\|_\mathbb{Y} \leq \lc f\rf \left(\|\varphi\|_\mathbb{X}\right).
\end{equation}
For any $r\in (0,+\infty)$, we let $\F_r(\mathbb{X},\mathbb{Y})$ be the space of formal series $f$
s.t. $\lc f \rf (r) < +\infty$. 
In the case where $\mathbb{Y} =\R$, we simply note $\F_r(\mathbb{X}):= \F_r(\mathbb{X},\mathbb{R})$

Lastly a family $(f_a)_{a\in A}$ of elements in $\mathbb{F}(\mathbb{X},\mathbb{Y})$ is called
a \emph{normal family of analytic maps of multiradius $r$} if there exists
$X\in \mathbb{F}(\R)$ s.t. $\rho(X)=r$ and,
setting $X(z) = \sum_{p=0}^\infty X^{(p)}z^p$,
$\forall a\in A$, $\forall p\in \N$, $\|f^{(p)}_a\|_\otimes \leq X^{(p)}$
(hence in particular $\rho_\otimes(f_a)\geq r$).\\

\noindent
\textbf{The nonlinearity} ---
We note $E_{(1)}:=  E\times E_2\times \mathcal{L}(\mathbb{R}^n,E)$.
We assume that the map $N$ is real analytic from $E_{(1)}$ to $E$
and that its multiradius of convergence is positive.

For applications to equations in Physics,
we are particularly interested in systems (\ref{equa0}) of the form\footnote{Actually any system of the form
$L_1u_1+N_1(u) = 0$ and $L_2u_2+\hat{N}_2(u,\partial u) = 0$ can be set in the form (\ref{L+N-a}) through
the substitution $N_2(u,\partial_0u_2,\vec{\partial}u):= \hat{N}_2(u,-N_1(u),\partial _0u_2,\vec{\partial}u)$.}:
\begin{equation}\label{L+N-a}
\left\{
\begin{array}{ccc}
 L_1u_1 + N_1(u) & = & 0\\
  L_2u_2 + N_2(u,\partial_0u_2,\vec{\partial}u) & = & 0
\end{array}
\right.
\end{equation}
where $N_1:E\longrightarrow E_1$ and $N_2:E_{(1)} \longrightarrow E_2$.
Motivated by the Yang--Mills system, we are led to consider the case where
\emph{$N_2$ is affine in $\partial u$}, i.e. there exist real analytic functions 
$J$ and $K_i^\mu:E\longrightarrow E_2$ s.t.
\begin{equation}\label{affine-hypothesis}
 N_2(u,\partial u) = J(u) + \sum_{i=1}^{\hbox{\footnotesize{dim}}E}\sum_{\mu=0}^nK_i^\mu(u)\partial_\mu u^i.
\end{equation}
By setting $N:= (N_1,N_2)$, we see that System (\ref{L+N-a}) is equivalent to (\ref{equa0}).

For any interval $I\subset \R$, we define
\begin{equation}\label{definition-Es}
 \mathcal{E}^s(I):=\{ u\in \mathcal{F}^s(I)|\
Lu + N(u,\partial u) = 0\hbox{ weakly}\}.
\end{equation}

\noindent
\textbf{The Lagrange--Duhamel vector field} ---
First define the `Green function' $G$ to be the unique distribution on $M$ with coefficients in $\hbox{End}(E)$,
which is a solution ot $LG=0$ and $L(YG) = \delta_0^{1+n}\otimes 1_E$,
where $Y(x):= 1_{[0,+\infty)}(x^0)$ is the Heaviside function.
Note that through the splitting $E = E_1\oplus E_2$, $G$ decomposes as
\begin{equation}\label{Gdecomposition}
 G = \left( \begin{array}{cc}
            G_1 & 0 \\
            0 & G_2
          \end{array} \right),
\end{equation}
where $G_1|_0 = \delta_0^n\otimes 1_{E_1}$, $G_2|_0 = 0$ and $\partial_0G_2|_0 = \delta_0^n\otimes 1_{E_2}$.

We define the time dependent  \emph{Lagrange--Duhamel} vector field $V_t$ on $\mathcal{E}_0^s(\R)$ by:
$\forall x\in M$,
\begin{equation}\label{lagrange-duhamel}
V_t(\varphi)(x):= \int_{\mathbb{R}^n} d\vec{y}\ G(x^0-t,\vec{x}-\vec{y})N(\varphi,\partial \varphi)(t,\vec{y})
= \int_{y^0=t}d\vec{y}\,G_y(x)N(\varphi,\partial\varphi)(y),
\end{equation}
where $G_y(x):= G(x-y)$.
By Theorem \ref{theo-magic-dynamics} a map $u$ is a solution of (\ref{equa0}) iff
\begin{equation}\label{magic-dynamics}
 \frac{d(\Theta_tu)}{dt} + V_t\left(\Theta_tu\right) = 0,
\end{equation}

\noindent
\textbf{The chronological exponential} --- The chronological exponential of $(V_t\cdot)_{t\in I}$
(if it exists) is the operator acting on $\mathbb{F}$ defined by
\begin{equation}\label{def-chronologic}
\begin{array}{ccl}
\displaystyle T\hbox{exp}\int_{t_1}^{t_2} d\tau V_\tau\cdot & := & \displaystyle 
\sum_{j=0}^\infty \int_{t_1<\tau_1<\cdots <\tau_j<t_2}d\tau_1\cdots d\tau_j(V_{\tau_j}\cdots V_{\tau_1}\cdot) ,
\quad \hbox{if }t_2>t_1\\
\hbox{or} & := & \displaystyle 
\sum_{j=0}^\infty \int_{t_2<\tau_j<\cdots <\tau_1<t_1}d\tau_1\cdots d\tau_j(-1)^j(V_{\tau_j}\cdots V_{\tau_1}\cdot) ,
\quad \hbox{if }t_2<t_1,
\end{array}
\end{equation}
with the convention that the first term in the sum ($j=0$) is the identity operator.

\section{The Lagrange--Duhamel vector field formulation}
The aim of this Section is to prove the following results.
\begin{lemm}\label{lemm-existence-Theta}
Let $J$ and $I$ be two intervals of $\R$ s.t. $J\subset I$.
Assume that $P_2$ satisfies (\ref{massive-nonmassive}) and that
either $\mu_0$ in (\ref{massive-nonmassive}) is positive or $I$ is bounded.
Then the map
\[
 \begin{array}{cccc}
  \Theta: & I\times \mathcal{F}^s(J) & \longrightarrow & \mathcal{E}^s_0(I)
 \end{array}
\]
defined by (\ref{composed-maps}) exists and is
continuous. Moreover there exists a constant $C_\Theta(I)>0$ s.t.
\[
 \forall u \in \mathcal{F}^s(J), \forall t\in J, \quad \|\Theta_tu\|_{\mathcal{F}^s(I)}
 \leq C_\Theta(I)\|u\|_{\mathcal{F}^s(J)}.
\]
\end{lemm}

\begin{prop}\label{proposition-fund-Section2}
Assume that $P_2$ satisfies (\ref{massive-nonmassive}) and that
either $\mu_0$ in (\ref{massive-nonmassive}) is positive or $I\subset \R$ is bounded. 
Then there exists some constant $Q_s>0$ such that the following holds.

Let $N:E_{(1)}\longrightarrow E$ be a real analytic map of multiradius of convergence $\rho_\otimes(N)>0$. Assume that: either
\begin{enumerate}
 \item $N_2$ does not depend on $\partial u$ and $s>n/2$; or
 \item $N_2$ is affine in $\partial u$, i.e. (\ref{affine-hypothesis}) holds and
 $s>n/2>s-r\geq 0$; or
 \item $s>n/2+1$.
\end{enumerate}
Then $\exists \rho_\otimes(V)$ s.t. $Q_s\rho_\otimes(V)\geq \rho_\otimes(N)$ and
$\forall (t,\varphi)\in I\times B_{\mathcal{E}_0^s(I)}(\rho_\otimes(V))$,
the quantity
\[
 V(t,\varphi) = V_t(\varphi):=
\int_{y^0=t} d\vec{y}\ G_yN(\varphi,\partial \varphi)(y)
\]
is well-defined. Moreover the map $V: I\times B_{\mathcal{E}_0^s}(\rho_\otimes(V))
\longrightarrow \mathcal{E}_0^s(I)$ is continuous and
$(V_t)_{t\in I}$ is a normal family of analytic maps of multiradius
equal to $\rho_{\otimes}(V)$.
\end{prop}

\begin{theo}\label{theo-magic-dynamics}
Assume the same hypotheses as in Proposition \ref{proposition-fund-Section2}.
Let $u\in \mathcal{F}^s(I)$ s.t. $\|u\|_{\mathcal{F}^s(I)}\leq \rho_{\otimes}(V)$.
Then $u$ belongs to $\mathcal{E}^s(I)$ (i.e. is a weak solution of $Lu + N(u,\partial u) = 0$) iff the map
\[
 \begin{array}{ccc}
  I & \longrightarrow & \mathcal{E}^s_0(I)\\
t & \longmapsto & \Theta_tu
 \end{array}
\]
is $\mathcal{C}^1$ and satisfies (\ref{magic-dynamics}), i.e.
$\frac{d(\Theta_tu)}{dt} + V_t\left(\Theta_tu\right) = 0$.
\end{theo}
\begin{rema}\label{remarksursegal1}
 Analogues of Proposition \ref{proposition-fund-Section2} and Theorem \ref{theo-magic-dynamics} can be
 proved without difficulty in the case where $E=E_2$ (i.e. $L$ is a fully second
 order operator), $1=s\leq n/2$ and if $N=N_2$ is a polynomial of degree
 $\hbox{deg}N\leq n/n-2$. This is a consequence of the Sobolev embedding
 $H^1(\mathbb{R}^n)\hookrightarrow L^{2n/n-2}(\mathbb{R}^n)$, which allows to estimate
 the nonlinearity in $H^0(\R^n) = L^2(\R^n)$. The relevant cases
 are: $n=2$ ($\hbox{deg}N$ is arbitrary); $n=3$ ($\hbox{deg}N \leq 3$) and
 $n=4$ ($\hbox{deg}N\leq 2$). The proof is left to the Reader. The special case
 $n=3$ and $N(u) = u^3$ will be treated in Section \ref{section-curved}.
\end{rema}

\subsection{Existence and continuity of $\Theta$}\label{paragraph2.2}
For any $(\psi,\chi)\in \hbox{Cau}^s$
and $t\in \R$, we recall that $\Phi_t(\psi,\chi)$ is equal to the unique solution $\varphi$ of $L\varphi=0$
on $\mathbb{R}^{n+1}$ s.t. $[\varphi]_t = (\psi,\chi)$. We also denote
by $\Phi_t(\psi,\chi)$ the restriction of this map to any subset $I\times \R^n$.
We set
\begin{equation}\label{big-composed-maps}
 \begin{array}{cccc}
  \Phi : & \R\times \hbox{Cau}^s & \longrightarrow & \mathcal{E}_0^s(I)
\\
 & (t,(\psi,\chi)) &  \longmapsto & \Phi_t(\psi,\chi).
 \end{array}
\end{equation}

\begin{prop}\label{proposition-continuity-of-Phi}
Assume that $P_2$ satisfies (\ref{massive-nonmassive}), then the linear map $\Phi$
defined  by (\ref{big-composed-maps}) is well-defined and continuous.
Assume that either $\mu_0$ in (\ref{massive-nonmassive}) is positive or $I\subset \R$ is bounded.
Then there exists a constant $C_\Phi(I)>0$ s.t.
\begin{equation}\label{estimate-on-Phi}
||\Phi_t(\psi,\chi)||_{\mathcal{F}^s(I)} \leq C_\Phi(I) ||(\psi,\chi)||_{\hbox{\tiny{Cau}}^s},
\quad
\forall (\psi,\chi)\in \Caus.
\end{equation}
\end{prop}
\emph{Proof} --- Since $\Phi$ is linear in its second argument we can decompose the problem
in two subcases and assume
either $L=L_1 = \gamma^0\partial_0 +P_1(\vec{\partial})$, or $L=L_2 = \partial_0^2 +P_2(\vec{\partial})$,
separately.

\noindent
\textbf{Case $L=L_1$:} we need to show that any solution $\varphi$ to
$\gamma^0\partial_0\varphi + P_1(\vec{\partial})\varphi = 0$, s.t. $\varphi|_t = \psi$
belongs to $\mathcal{F}^s(I)$ and depends
continuously on $(t,\psi)$, where $\psi\in H^s(\mathbb{R}^n,E)$.
Setting $\epsilon(\xi):= i(\gamma^0)^{-1}P_1(i\xi)$, the equation reads
$\partial_0 \widehat{\varphi}(\tau,\xi) - i\epsilon(\xi)\widehat{\varphi}(\tau,\xi) = 0$,
its solution is given by $\widehat{\varphi}(\tau,\xi) = e^{i\epsilon(\xi)(\tau-t)}\widehat{\psi}(\xi)$ and
its time derivative by
$\partial_0\widehat{\varphi}(\tau,\xi) = i\epsilon(\xi)e^{i\epsilon(\xi)(\tau-t)}\widehat{\psi}(\xi)$.
The result then follows by standard majorations and Lebesgue's dominated theorem.

\noindent
\textbf{Case $L=L_2$:} We need to show that the solution $\varphi$ of $(\partial_0)^2\varphi + P_2(\vec{\partial})\varphi =0$,
s.t. $\varphi|_t = \psi$ and $\partial_0\varphi|_t = \chi$ depends continuously on
$(t,\psi,\chi)\in \R\times H^s(\mathbb{R}^n,E)\times H^{s-r}(\mathbb{R}^n,E)$.
Assuming that $P_2(i\xi)$ is positive Hermitian
and setting $\epsilon(\xi):= \sqrt{P_2(i\xi)}$, the equation reads
$(\partial_0)^2 \widehat{\varphi}(\tau,\xi)  + \epsilon(\xi)^2\widehat{\varphi}(\tau,\xi) = 0$.
Its solution is $\widehat{\varphi}(\tau,\xi) = \cos \epsilon(\xi)(\tau-t)\widehat{\psi}(\xi)
+ \epsilon(\xi)^{-1}\sin\epsilon(\xi)(\tau-t)\chi(\xi)$ and its time derivative is
$\partial_0\widehat{\varphi}(\tau,\xi) = -\epsilon(\xi)\sin \epsilon(\xi)(\tau-t)\widehat{\psi}(\xi)
+ \cos\epsilon(\xi)(\tau-t)\chi(\xi)$.
The proof that $\varphi\in \mathcal{E}_0^s(I)$ and its continuous dependence on $(t,\psi,\chi)$ follows the same
lines as for first order equations. However the factor $\epsilon(\xi)^{-1}$ in the expression of $\widehat{\varphi}$
may pose a slight difficulty in proving that $\varphi$
is in $\mathcal{C}^0(I,H^s(\mathbb{R}^n,E))$ and that it depends continuously in $t$.
In the `massive case' (i.e. ${\mu_0}$ in (\ref{massive-nonmassive}) is positive) this difficulty
does not occur because of the inequality $|\epsilon(\xi)^{-1}|\leq (\alpha{\mu_0})^{-1/2}$.
In the `non massive' case (i.e. ${\mu_0}$ in (\ref{massive-nonmassive}) vanishes)
we only have $|\epsilon(\xi)^{-1}|\leq \alpha^{-1/2}|\xi|^{-r}$. However by using the inequality:
\begin{equation}\label{A-ptit-lemme}
 \left|\frac{\sin\epsilon(\xi) t}{\epsilon(\xi)}\right| \leq \frac{\sqrt{1+t^2}}{\sqrt{1+ \alpha|\xi|^{2r}}},
\end{equation}
we can prove the result by working with $\mathcal{E}^s_0$ endowed with the norm
\[
 \|u\|_{\check{L}^\infty \caus}:= \|u\|_{\check{L}^\infty(I, H^s)}
+ \|\partial_0u_2\|_{L^\infty (I,H^{s-r})},
\]
where $\|u\|_{\check{L}^\infty(I, H^s)}:= \sup_{\tau\in I}\frac{\|u|_\tau\|_{H^s}}{\sqrt{1+\tau^2}}$.
The conclusion follows if $I$ is bounded, since then both norms $L^\infty\Caus$
and $\check{L}^\infty \Caus$ are equivalent.  \hfill $\square$\\

\noindent
\emph{Proof of Lemma \ref{lemm-existence-Theta}} ---
Since the map $\Theta$ is obtained by composing $\Phi$ with the map
$I\times \mathcal{F}^s(J)\longrightarrow I\times \hbox{Cau}^s$,
$(t,u)\longmapsto (t,[u]_t)$, which is obviously continuous, Lemma \ref{lemm-existence-Theta}
is a straightforward consequence of Proposition \ref{proposition-continuity-of-Phi}.\hfill $\square$

\subsection{Estimate on the nonlinearity}

\emph{The goal of this section is to collect results to prove Proposition \ref{proposition-fund-Section2}.}
As a preliminary result we prove 

\begin{lemm}\label{lemma-series-are-continuous}
Let $\mathbb{X}$ and $\mathbb{Y}$ be Banach spaces.
Let $f= \sum_{p=0}^\infty f^{(p)}$ be a formal series from $\mathbb{X}$ to $\mathbb{Y}$.
Assume that its multiradius of convergence $\rho_\otimes(f) > 0$.
Then the map $f$ defined by $\forall \varphi\in \mathbb{X}$,
$f(\varphi) = \sum_{p=0}^\infty f^{(p)}(\varphi)$ is $\mathcal{C}^\infty$ on
$B_\mathbb{X}(\rho_\otimes(f))$.

In particular:
$\forall r$ s.t. $0<r<\rho_\otimes(f)$, $\forall \varphi,\psi\in \mathbb{X}$, $\forall h\in \R$ s.t.
$\|\varphi\|_\mathbb{X},\|\varphi + h\psi\|_\mathbb{X}\leq r$
\begin{equation}\label{lipschitzlemma}
 \|f(\varphi + h\psi) - f(\varphi)\|_\mathbb{Y}
 \leq \frac{d\lc f\rf}{dz}(r)\|h\psi\|_\mathbb{X}.
\end{equation}
and
\begin{equation}\label{seriesC2}
 \left\|f(\varphi+h\psi)-f(\varphi) - h\delta f_\varphi(\psi)\right\|_{\mathbb{Y}}
 \leq \frac{1}{2} \frac{d^2\lc f\rf}{dz^2}(r)h^2\|\psi\|_{\mathbb{X}}
\end{equation}

\end{lemm}
\emph{Proof} --- We prove only that $f$ is $\mathcal{C}^1$ and hence (\ref{lipschitzlemma}) and
(\ref{seriesC2}) and leave the general case to the Reader. We first prove (\ref{lipschitzlemma}).
Let $p\in \N$. From the identity
$f^{(p)}(\varphi + h\psi) - f^{(p)}(\varphi) =
\sum_{j=1}^p  f^{(p)}_{\otimes}((\varphi + h\psi)^{\otimes j-1}\otimes
h\psi\otimes \varphi^{\otimes p-j})$
we deduce
\[
 \|f^{(p)}(\varphi + h\psi) - f^{(p)}(\varphi)\|_\mathbb{Y}\leq
 \sum_{j=1}^p \|f^{(p)}\|_{\otimes}
 \|\varphi + h\psi\|_\mathbb{X}^{j-1} \| h\psi \|_\mathbb{X}
\|\varphi\|_\mathbb{X}^{p-j}.
\]
Thus if $\|\varphi\|_\mathbb{X},\|\varphi + h\psi\|_\mathbb{X}\leq r$,
\begin{equation}\label{F-continu}
\|f^{(p)}(\varphi + h\psi) - f^{(p)}(\varphi)\|_\mathbb{Y}
 \leq p\|f^{(p)}\|_{\otimes}r^{p-1}\|h\psi \|_\mathbb{X}.
\end{equation}
Hence by summing up on $p\in\N$ and using
$\frac{d\lc f\rf}{dz}:= \sum_{p=0}^\infty p\|f^{(p)}\|_{\otimes}z^{p-1}$, we deduce that 
(\ref{lipschitzlemma}) holds if $r<\rho_\otimes(f)$.

The proof of (\ref{seriesC2}) is similar. We start from the identity
\[
 f^{(p)}(\varphi+h\psi)-f^{(p)}(\varphi) - p f^{(p)}_\otimes(h\psi\otimes \varphi^{\otimes p-1})
 = \sum_{j_1=1}^p\sum_{j_2=1}^{j_1-1}f^{(p)}_\otimes \left((\varphi+h\psi)^{\otimes j_2-1}
 \otimes (h\psi)^{\otimes 2}\otimes \varphi^{\otimes p-j_2-1}\right),
\]
%
%
from which we deduce that, if $\|\varphi\|_\mathbb{X},\|\varphi + h\psi\|_\mathbb{X}\leq r$,
\[
 \left\|f^{(p)}(\varphi+h\psi)-f^{(p)}(\varphi) - p f^{(p)}_\otimes(h\psi\otimes \varphi^{\otimes p-1})\right\|_{\mathbb{Y}}
 \leq \frac{p(p-1)}{2}\|f^{(p)}\|_\otimes r^{p-2}\|h\psi\|_\mathbb{X}^2.
\]
Hence (\ref{seriesC2}) follows by summing up on $p\in \N$.
\hfill $\square$\\

Let $V$ and $W$ be two real vector spaces of (finite) dimension $d_V$ and $d_W$ respectively, let $k\in \N$ and 
$L\in \mathcal{Q}^k(V,W)$. Let $L_\otimes:V^{\otimes k}\longrightarrow W$
the associated polarized linear map. Using bases on $V$ and $W$, $L_{\otimes}$ has
the coordinates representation:
\begin{equation}\label{def-L-polarized}
 L_{\otimes}^i(z_1\otimes \cdots\otimes z_k):= \sum_{j_1,\cdots,j_k=1}^{d_V}
L_{j_1,\cdots,j_k}^iz_1^{j_1}\cdots z_k^{j_k}
\quad \forall i=1,\cdots, d_W, \forall z_1,\cdots ,z_k \in V.
\end{equation}
where, $\forall a$ s.t. $1\leq a\leq k$, $\left(z^j_a\right)_{1\leq j\leq d_V}$ are the coordinates of $z_a$
and the coefficients $L_{j_1,\cdots,j_k}^i$ are symmetric in $(j_1,\cdots,j_k)$.
We set
\begin{equation}\label{def-L-norm}
 |L|:= \sum_{i=1}^{d_W}\sup_{1\leq j_1,\cdots,j_k\leq d_V}\left|L_{j_1,\cdots,j_k}^i\right|.
\end{equation}
One can easily check that (see (\ref{def-tensornorm}))
\begin{equation}\label{encadrement-normes-de-L}
 \|L\|_\otimes \leq |L| \leq d_W\|L\|_\otimes.
\end{equation}
The following result uses the fact that, if $s>n/2$,
then $H^s(\mathbb{R}^n)$ is an algebra, i.e. the product of two functions $f,g\in H^s(\mathbb{R}^n)$
belongs to $H^s(\mathbb{R}^n)$ and there exists a constant $Q_s $ s.t.
$||f||_{H^s}||g||_{H^s} \leq Q_s  ||fg||_{H^s}$.
\begin{lemm}\label{lemmaL-preN1}
Let $V$ and $W$ be two real vector spaces of finite dimension, let $k\in \N$ and 
let $L\in \mathcal{Q}^k(V,W)$. Assume that $s>n/2$.
Then one can define the homogeneous polynomial map $\mathcal{L}\in \mathcal{Q}^k\left((H^s(\mathbb{R}^n,V))^{\otimes k},
H^s(\mathbb{R}^n,W)\right)$ by $\mathcal{L}_\otimes (\varphi_1\otimes\cdots \otimes\varphi_k)(x):=
L_\otimes (\varphi_1(x)\otimes\cdots \otimes\varphi_k(x))$ a.e., this map is linear continuous and
satisfies
\begin{equation}\label{estimate-norm-of-L}
\|\mathcal{L}\|^{\otimes}\leq Q_s^{k-1}|L|.
\end{equation}
\end{lemm}
\emph{Proof} --- A straightforward recursion shows that
$\|\varphi_1^{j_1}\cdots \varphi_k^{j_k}\|_{H^s}\leq
Q_s^{k-1}\|\varphi_1^{j_1}\|_{H^s}\cdots \|\varphi_k^{j_k}\|_{H^s}$ and
hence:
\[
 \begin{array}{ccl}
  \displaystyle \|\mathcal{L}(\varphi_1\otimes\cdots \otimes\varphi_k)\|_{H^s} & = &
\displaystyle \sum_{i=1}^{d_W} \|\mathcal{L}^i(\varphi_1\otimes\cdots \otimes\varphi_k)\|_{H^s}\\
& \leq &
\displaystyle \sum_{i=1}^{d_W}\sum_{j_1,\cdots,j_k=1}^{d_V}
\left|L_{j_1,\cdots,j_k}^i\right| \|\varphi_1^{j_1}\cdots \varphi_k^{j_k}\|_{H^s}\\
& \leq &
\displaystyle \sum_{i=1}^{d_W}\sup_{1\leq j_1,\cdots,j_k\leq d}
\left|L_{j_1,\cdots,j_k}^i\right|\sum_{j_1,\cdots,j_k=1}^{d_V}
Q_s^{k-1}\|\varphi_1^{j_1}\|_{H^s}\cdots \|\varphi_k^{j_k}\|_{H^s}.
 \end{array}
\]
Hence the result follows by using (\ref{def-L-norm}).\hfill $\square$\\

In the following we use the notations:
\begin{equation}\label{mapCaustoHsE1}
 \begin{array}{ccc}
  \hbox{Cau}^s & \longrightarrow & H^s(\R^n,E) \times H^{s-r}(\R^n,E_2) \times H^{s-1}(\R^n,\mathcal{L}(\R^n,E))\\
  \phi= (\varphi,\chi) & \longmapsto & \phi^{(1)} = (\varphi,\chi,\partial_1\varphi,\cdots ,\partial_n\varphi).
 \end{array}
\end{equation}
As a first application of Lemma \ref{lemmaL-preN1}, given $N_1=\sum_{k=0}^\infty N_1^{(k)}\in \F(E,E_1)$,
we define for any $k\in \N$ the map $(\mathcal{N}_1^{(k)})_\otimes: (\hbox{Cau}^s)^{\otimes k} \longrightarrow  H^s(\R^n,E_1)$
by $\forall (\phi_1,\cdots ,\phi_k)\in (\hbox{Cau}^s)^k$,
\begin{equation}\label{definitioncalligraphicN1}
\mathcal{N}_1^{(k)}(\phi_1\otimes\cdots \otimes\phi_k)(x):=
N^{(k)}_1(\varphi_1(x)\otimes \cdots \otimes \varphi_k(x)),\quad \hbox{for a.e. }x\in \R^n.
\end{equation}
We then deduce from (\ref{estimate-norm-of-L}) the estimate
$\|\mathcal{N}_1^{(k)}\|^{\otimes}\leq Q_s^{k-1}|N_1^{(k)}|$. A similar estimate can be obtained
for $N_2$ if this function does not depend on $\partial u$ or if $s>n/2+1$.

However if $N_2$ is affine in $\partial u$, i.e. has the form (\ref{affine-hypothesis}) and if we suppose that
$0\leq s-r<n/2<s$, then we use the fact that the product $(f,g)\longmapsto fg$ also maps continuously
$H^s(\mathbb{R}^n) \times H^{s-r}(\mathbb{R}^n)$ to $H^{s-r}(\mathbb{R}^n)$ and that there
exists a constant $q_{s,r}$ s.t.
\begin{equation}\label{estimate-paraproduct}
\|fg\|_{H^{s-r}}\leq q_{s,r}Q_s\|f\|_{H^s}\|g\|_{H^{s-r}}.
\end{equation}
This can be proved by splitting the product $fg$ as
the sum $T_fg + T_gf + R(f,g)$, where $(f,g)\longmapsto T_fg$ is the paraproduct and by
estimating each term separately:
$T_fg\in H^{s-r}(\mathbb{R}^n)$ because $f\in H^s(\mathbb{R}^n)\subset L^\infty(\mathbb{R}^n)$,
$T_gf\in H^{s+(s-r)-n/2}(\mathbb{R}^n) \subset H^{s-r}(\mathbb{R}^n)$ because $s-r<n/2$ and
$R(f,g)\in H^{s+(s-r)-n/2}(\mathbb{R}^n) \subset H^{s-r}(\mathbb{R}^n)$ because $s+(s-r)>n/2$ (see \cite{alinhac-gerard},
Exercise A.5, page 109).

For the following remind the notation introduced in (\ref{mapCaustoHsE1}).
We also use the notation $\delta_0\phi:= \chi$, $\delta_i\phi:= \partial_i\varphi$, for $1\leq i\leq n$,
$\delta\phi:= (\delta_\mu\phi)_{0\leq \mu \leq n}$, $\forall \phi= (\varphi,\chi)\in \hbox{Cau}^s$.
\begin{lemm}\label{lemmaN2}
Let $N_2^{(k)}\in \mathcal{Q}^k(E_{(1)}, E_2)$ satisfying (\ref{affine-hypothesis})
with $J^{(k)}\in \mathcal{Q}^k(E, E_2)$ and $K_i^{(k)\mu}\in \mathcal{Q}^{k-1}(E_{(1)}, E_2)$.
Assume that $0\leq s-r<n/2<s$.
Then one can define the map $\mathcal{N}_{2\otimes}^{(k)}$ from $(\hbox{Cau}^s)^{\otimes k}$ to
$H^{s-r}(\R^n,E_2)$ by
$\mathcal{N}_2^{(k)}\left(\phi_1\otimes \cdots \otimes \phi_k\right)(x):=
N_2^{(k)}\left(\phi_1^{(1)}(x)\otimes\cdots \otimes\phi_k^{(1)}(x)\right)$, for a.e. $x\in \R^n$, this map
is linear continuous and satisfies:
\begin{equation}\label{estimate-N2-multi}
 \|\mathcal{N}_2^{(k)}\|_{\otimes} \leq  Q_s^{k-1}\left(|J^{(k)}|+q_{r,s}\sqrt{n}|K^{(k)}|\right)
\quad \hbox{for}\quad
|K^{(k)}|:= \sup_{1\leq i\leq \hbox{\footnotesize{dim}}E}\sup_{0\leq \mu\leq n}|K^{(k)\mu}_i|.
\end{equation}
\end{lemm}
\emph{Proof} --- We start from the decomposition
\begin{equation}\label{starting-decomp}
\begin{array}{ccl}
\displaystyle \mathcal{N}_{2}^{(k)}\left(\phi_1\otimes\cdots \otimes\phi_k\right)
& = & \displaystyle \mathcal{J}^{(k)}\left(\varphi_1\otimes\cdots \otimes\varphi_k\right)\\
& & \displaystyle + \frac{1}{k}\sum_{a=1}^k\sum_{i=1}^{\hbox{\footnotesize{dim}}E}\sum_{\mu=0}^n
\mathcal{K}^{(k)\mu}_{i}\left(\varphi_1\otimes\cdots \otimes\widehat{\varphi}_a\otimes\cdots \otimes\varphi_k\right)\delta_\mu\phi^i_a.
\end{array}
\end{equation}
Inequality (\ref{estimate-norm-of-L}) gives us $\|\mathcal{J}^{(k)}\left(\varphi_1\otimes\cdots \otimes\varphi_k\right)\|_{H^s}\leq 
Q_s^{k-1}|J^{(k)}| \|\varphi_1\|_{H^s}\cdots \|\varphi_k\|_{H^s}$. This implies automatically a similar estimate
in $H^{s-r}$. The $H^{s-r}$ norm of the r.h.s. term in (\ref{starting-decomp}) is estimated by
using (\ref{estimate-paraproduct}):
\[
 \begin{array}{cl}
  \leq & \displaystyle \frac{q_{r,s}Q_s}{k}\sum_{a=1}^k
\left(\sup_{1\leq i\leq \hbox{\footnotesize{dim}}E}\sup_{0\leq \mu\leq n}
\left\|\mathcal{K}^{(k)\mu}_{i}\left(\varphi_1\otimes\cdots \otimes\widehat{\varphi}_a\otimes\cdots \otimes\varphi_k\right)\right\|_{H^s}\right)
\left(\sum_{i=1}^{\hbox{\footnotesize{dim}}E}\sum_{\mu=0}^n\|\delta_\mu\phi^i_a\|_{H^{s-r}}\right)\\
\leq & \displaystyle \frac{q_{r,s}Q_s}{k}\sum_{a=1}^k
\left(|K^{(k)}|Q_s^{k-2}\|\varphi_1\|_{H^s}\cdots\widehat{\|\varphi_a\|_{H^s}}\cdots\|\varphi_k\|_{H^s}\right)
\left(\sum_{\mu=0}^n\|\delta_\mu\phi_a\|_{H^{s-r}}\right),
 \end{array}
\]
where we have used (\ref{estimate-norm-of-L}). However by using Cauchy--Schwarz inequality and the fact
that $r\geq 1$ we have $\forall a = 1,\cdots ,k$, 
\[
 \begin{array}{ccl}
\displaystyle  \sum_{\mu=0}^n\|\delta_\mu\phi_a\|_{H^{s-r}} & \leq &\displaystyle  \|\chi_a\|_{H^{s-r}}
+ \sqrt{n}\left(\sum_{\mu=0}^n\|\partial_\mu\varphi_a\|_{H^{s-r}}^2\right)^{1/2}\\
& \leq  & \displaystyle \|\chi_a\|_{H^{s-r}} + \sqrt{n}\|\varphi_a\|_{H^{s+1-r}}
\leq \sqrt{n}\|\phi_a\|_{\caus}.
 \end{array}
\]
Hence the $H^{s-r}$ norm of the r.h.s. term in (\ref{starting-decomp}) is estimated by:
\[
 \begin{array}{cl}
\leq & \displaystyle Q_s^{k-1}\frac{q_{r,s}\sqrt{n}}{k}|K^{(k)}|\sum_{a=1}^k
\left(\|\varphi_1\|_{H^s}\cdots\widehat{\|\varphi_a\|_{H^s}}\cdots\|\varphi_k\|_{H^s}\right)
\|\phi_a\|_{\caus}\\
\leq & \displaystyle Q_s^{k-1}q_{r,s}\sqrt{n}|K^{(k)}|\|\phi_1\|_{\caus}\cdots \|\phi_k\|_{\caus}.
 \end{array}
\]
Hence (\ref{estimate-N2-multi}) follows.\hfill $\square$\\

\noindent
Let's summarize Lemmas \ref{lemmaL-preN1} and \ref{lemmaN2}. We can define for any $k\in \N$ the map
$\mathcal{N}^{(k)}_\otimes:(\hbox{Cau}^s)^{\otimes k}\longrightarrow \hbox{Cau}^s$ by:
for a.e. $x\in \R^n$,
\[
 \mathcal{N}^{(k)}\left(\phi_1\otimes \cdots \otimes \phi_k\right)(x):=
 \left(\iota\circ \mathcal{N}^{(k)}_1\left(\phi_1\otimes \cdots \otimes \phi_k\right)(x),
 \mathcal{N}^{(k)}_2\left(\phi_1\otimes \cdots \otimes \phi_k\right)(x)\right),
\]
where $\iota: E_1\longrightarrow E$ is the natural inclusion.
Remark that\\
$\|\mathcal{N}^{(k)}(\phi_1\otimes\cdots\otimes \phi_k)\|_{\caus}
= \|\mathcal{N}_1^{(k)}(\varphi_1\otimes\cdots\otimes \varphi_k)\|_{H^s}
        + \|\mathcal{N}_2^{(k)}(\phi_1\otimes\cdots\otimes \phi_k)\|_{H^{s-r}}$.
\begin{prop}\label{prop-N1+N2}
Assume that $N=\sum_{k=0}^\infty N^{(k)}\in \F(E_{(1)},E)$ satisfies the same
hypotheses as in Proposition \ref{proposition-fund-Section2}.
Then, for any $k\in \N$, the linear map $\mathcal{N}^{(k)}_\otimes$ from $(\hbox{Cau}^s)^{\otimes k}$
to $\hbox{Cau}^s$ is continuous and satisfies :
\begin{equation}\label{lien-entre-rayon}
 \|\mathcal{N}^{(k)}\|_{\otimes}\leq Q_s^kC(s,r,n)\|N^{(k)}\|_{\otimes}.
\end{equation}
Hence if $\rho_\otimes(N)>0$, then $\mathcal{N} =
\sum_{k=0}^\infty \mathcal{N}^{(k)}\in \mathbb{F}_{Q_s^{-1}\rho_\otimes(N)}(\Caus,\Caus)$.
\end{prop}
\emph{Proof} --- Case (i) where $s>n/2$ and $N_2^{(k)}$ does not depend on $\partial u$
and Case (iii) where $s>n/2+1$ are similar and can dealt by
applying Lemma \ref{lemmaL-preN1} for $L=N_1^{(k)}$ and $L=N_2^{(k)}$. We then obtain
$\| \mathcal{N}^{(k)}\|_{\otimes}\leq Q_s^{k-1}(|N_1^{(k)}| + |N_2^{(k)}|)$.
In Case (ii), we apply Lemma \ref{lemmaL-preN1} for $L=N_1^{(k)}$ and Lemma \ref{lemmaN2} for $N_2^{(k)}$ to get
$\| \mathcal{N}^{(k)}\|_{\otimes}\leq Q_s^{k-1}(|N_1^{(k)}| + |J^{(k)}| + q_{r,s}\sqrt{n}|K^{(k)}|)$.
In any case (\ref{lien-entre-rayon}) follows by applying (\ref{encadrement-normes-de-L}) to
$L = N_1^{(k)}, N_2^{(k)}, J^{(k)},K^{(k-1)}$.
As a consequence $\sum_{k=0}^\infty \|\mathcal{N}^{(k)}\|_\otimes R^k$ converges if
$Q_sR<\rho_\otimes(N)$. \hfill $\square$\\

A first consequence of Proposition \ref{prop-N1+N2} is:
\begin{prop}\label{proposition-continuity-N(u)}
Let $s\in \R$ and $u\in \mathcal{F}^s(I)$.
Assume that $Q_s||u||_{\mathcal{F}^s(I)} < \rho_\otimes(N)$.
Assume that $N$ satisfies the same hypotheses in Proposition \ref{proposition-fund-Section2}.
Then $ \mathcal{N}(u,\partial u)\in \mathcal{C}^0(I,\hbox{Cau}^s)$, i.e. 
$ \mathcal{N}_1(u)\in \mathcal{C}^0(I,H^s(\mathbb{R}^n,E))$ and
$ \mathcal{N}_2(u,\partial u)\in \mathcal{C}^0(I,H^{s-r}(\mathbb{R}^n,E_2))$. 
\end{prop}
\emph{Proof} --- A straightforward consequence of Lemma \ref{lemma-series-are-continuous},
Proposition \ref{prop-N1+N2} and the continuity of $t\longmapsto [u]_t\in \hbox{Cau}^s$.\hfill $\square$\\

\subsection{The Lagrange--Duhamel vector field}
\emph{We prove Proposition \ref{proposition-fund-Section2}} --- For any interval $I\subset \R$ and
$k\in \N$, we define
$V_{\otimes}^{(k)}: I\times (\mathcal{E}_0^s(I))^{\otimes k} \longrightarrow  \mathcal{E}_0^s(I)$
by
\[
 V^{(k)}(t,\varphi_1\otimes \cdots \otimes \varphi_k) = V^{(k)}_t(\varphi_1\otimes \cdots \otimes \varphi_k)
:= \Phi_t\left( \mathcal{N}^{(k)}([\varphi_1]_t\otimes\cdots\otimes [\varphi_k]_t)\right)
\]
$V_\otimes^{(k)}$ is continuous since it is the composition of the maps
$I\times (\mathcal{E}^s_0(I))^{\otimes k}\ni (t,\varphi_1\otimes\cdots\otimes\varphi_k)
\longmapsto (t,[\varphi_1]_t\otimes\cdots\otimes [\varphi_k]_t)\in I\times (\Caus)^{\otimes k}$
(see (\ref{trivialmap})), $\mathcal{N}_\otimes^{(k)}$ (see Proposition \ref{prop-N1+N2})
and $\Phi$ (see Proposition \ref{proposition-continuity-of-Phi}). By using
(\ref{estimate-on-Phi}) and (\ref{lien-entre-rayon}) we deduce:
\begin{equation}\label{estimatefornormal}
 \|V^{(k)}\|_{\otimes}:=
 \sup_{t\in I} \|V^{(k)}_t\|_{\otimes}\leq
 C_\Phi(I)\|\mathcal{N}^{(k)}\|_{\otimes}
\leq C_\Phi(I) C(s,r,n)Q_s^k \|N^{(k)}\|_{\otimes}.
\end{equation}
Setting $V_t:= \sum_{k=0}^\infty V_t^{(k)}$,
(\ref{estimatefornormal}) implies that $(V_t)_{t\in I}$ is a normal family of analytic maps of multiradius 
of convergence $\rho_\otimes(V) \geq \rho_\otimes(N)/Q_s$.

To prove that $V:= \sum_{k=0}^\infty V^{(k)}$ is continuous on $I\times B_{\mathcal{E}^s_0(I)}(0,\rho_{\otimes}(V))$,
let $t,\tilde{t}\in I$, $\tilde{\varphi},\varphi\in B_{\mathcal{E}^s_0(I)}(r)$, where
$r<\rho_{\otimes}(V)$ and let us start from the inequality
\begin{equation}\label{vt2tildevarphimoinsvt1}
 \|V(\tilde{t},\tilde{\varphi}) - V(t,\varphi)\|_{\mathcal{E}^s_0(I)} \leq
 \| \Phi_{\tilde{t}}\left(\mathcal{N}([\tilde{\varphi}]_{\tilde{t}}) - \mathcal{N}([\varphi]_{t})\right)\|_{\mathcal{E}^s_0(I)}
  + \|(\Phi_{\tilde{t}} - \Phi_{t})\left(\mathcal{N}([\varphi]_{t})\right)\|_{\mathcal{E}^s_0(I)}.
\end{equation}
Fix $t,\varphi$ and $\varepsilon>0$, then we deduce from Proposition \ref{proposition-continuity-of-Phi}
that, by choosing $\tilde{t}$ sufficiently close to $t$, the last term in the r.h.s. of
(\ref{vt2tildevarphimoinsvt1}) is less than $\varepsilon/2$.

The first term in the r.h.s. of (\ref{vt2tildevarphimoinsvt1}) can be estimated
by using first (\ref{estimate-on-Phi}) and second (\ref{lipschitzlemma}):
\[
\begin{array}{ccl}
 \|\Phi_{\tilde{t}}(\mathcal{N}([\tilde{\varphi}]_{\tilde{t}}) - \mathcal{N}([\varphi]_{t}))\|_{\mathcal{E}^s_0(I)}
& \leq &
 C_\Phi(I)\|\mathcal{N}([\tilde{\varphi}]_{\tilde{t}}) - \mathcal{N}([\varphi]_{t})\|_{\caus}\\
 & \leq & C_\Phi(I) \frac{d\lc \mathcal{N}\rf}{dz}(r) \|[\tilde{\varphi}]_{\tilde{t}} - [\varphi]_{t}\|_{\caus},
\end{array}
\]
and hence will also be smaller than $\varepsilon/2$ if we choose $|\tilde{t}-t|$ and
$\|\tilde{\varphi}-\varphi\|_{\caus}$ sufficiently small.

To conclude observe that
\[
V_t^{(k)}(\varphi_1\otimes\cdots\otimes\varphi_k)(x)
 = \int_{y^0=t}d\vec{y}G(x-y)N^{(k)}((\varphi_1,\partial \varphi_1)\otimes\cdots\otimes
 (\varphi_k,\partial \varphi_k))(y).
\]
This can be proven by using the properties of $G$ (see (\ref{Gdecomposition})).
Hence in particular
\begin{equation}\label{lagrange-duhamel_k}
V_t^{(k)}(\varphi)(x):= \int_{\mathbb{R}^n} d\vec{y}
\ G(x^0-t,\vec{x}-\vec{y})N^{(k)}(\varphi,\partial \varphi)(t,\vec{y}).
\end{equation}
Formula (\ref{lagrange-duhamel}) which is similar to (\ref{lagrange-duhamel_k}) follows straightforwardly.
\hfill $\square$


\subsection{Derivability of $\Theta$}
We recall below Duhamel's formula (\ref{duhamelforever}).
Recall that $G$ is the distribution defined in (\ref{Gdecomposition}).
A generalization of Duhamel's formula for $L=L_2=\square_g +m^2$ on a curved pseudo Riemannian manifold will
given and proved in Section \ref{section-curved}.
\begin{prop}\label{proposition-duhamel}
Let $f=(f_1,f_2)$ where $f_1\in L_{loc}^1(\R, H^s(\R^n,E_1))$ and
$f_2\in L_{loc}^1(\R, H^{s-r}(\R^n,E_2))$.
Assume that $u\in \mathcal{F}^s(I)$ is a solution of $Lu = f$. Then 
\begin{equation}\label{duhamelforever}
\forall x\in \mathbb{R}^{1+n},\quad
 u(x) = \Theta_tu(x) + \int_t^{x^0} dy^0\int_{\mathbb{R}^n}d\vec{y}\ G(x-y)f(y).
\end{equation}
\end{prop}
We are now in position to give the:\\
\emph{Proof of Theorem \ref{theo-magic-dynamics}} --- The key observation is that $[\Theta_tu]_t = [u]_t$
implies $N\left(\Theta_tu,\partial (\Theta_tu)\right)|_t = N(u,\partial u)|_t$ and thus
\begin{equation}\label{simple-remark}
V_t(\Theta_tu)(x):= \int_{\mathbb{R}^n} d\vec{y}\ G(x^0-t,\vec{x}-\vec{y})N(u,\partial u)(t,\vec{y})
\end{equation}
Now since $u$ is a solution of $Lu+N(u,\partial u)=0$, we deduce from Proposition \ref{proposition-duhamel}
that $u(x) = \Theta_tu(x) - \int_t^{x^0} dy^0\int_{\mathbb{R}^n}d\vec{y}\ G(x-y)N(u,\partial u)(y)$,
which gives us thank to (\ref{simple-remark}):
\[
\forall x\in I\times \mathbb{R}^n,\quad u(x) = \Theta_tu(x) - \int_t^{x^0} dy^0V_{y^0}(\Theta_{y^0}u)(x).
\]
This implies the following identity $\forall t_1,t_2\in I$:
\begin{equation}\label{formule2.6}
 \Theta_{t_2}u - \Theta_{t_1}u + \int_{t_1}^{t_2} dy^0V_{y^0}(\Theta_{y^0}u) = 0
\quad \hbox{in }\mathcal{E}_0^s.
\end{equation}
Lastly Lemma \ref{lemm-existence-Theta} and Proposition \ref{proposition-fund-Section2} imply that
$I\ni t \longmapsto V_t(\Theta_tu) \in \mathcal{E}^s_0(I)$ is continuous. Hence
(\ref{formule2.6}) implies that $I\ni t \longmapsto \Theta_tu \in \mathcal{E}^s_0(I)$ is $\mathcal{C}^1$ and satisfies
(\ref{magic-dynamics}).
\hfill $\square$

\section{Curved space-times}\label{section-curved}
We show here how Theorem \ref{theo-magic-dynamics} can be extended to field equations on a curved space-time.
Let $\mathcal{M}$ be smooth manifold equipped with a $\mathcal{C}^\infty$ pseudo-Riemannian
metric $g$ of signature $(+,-,\cdots,-)$.
We denote by $\square_g = |g|^{-1}\partial_\mu\left(|g|g^{\mu\nu}\partial_\nu \right)$,
where $|g|:=  \sqrt{|\hbox{det}(g_{\mu\nu})|}$, the wave
operator and set  $L_g:= \square_g + m^2$.
A frame $(e_0,\dots,e_n)$ is called \emph{$g$-orthonormal} if $\langle e_\mu,e_\nu\rangle_g = \eta_{\mu\nu}$,
where $\eta_{00}=1$, $\eta_{ii}=-1$ if $i\neq 0$ and $\eta_{\mu\nu}=0$ if $\mu\neq \nu$.
We consider the non homogeneous scalar wave (or Klein--Gordon) equation on $\mathcal{M}$:
\begin{equation}\label{KG-on-curved}
L_gu:= |g|^{-1}\partial_\mu\left(|g|g^{\mu\nu}\partial_\nu u\right) + m^2u = f.
\end{equation}
Homogeneous (i.e. for $f=0$) solutions $u$ to (\ref{KG-on-curved}) are the critical points of the action functional
\[
 \mathcal{A}(u) = \frac{1}{2}\int_\mathcal{M}\left[|\partial u|_g^2 + m^2u^2\right]d\hbox{vol}_g,
\]
where $d\hbox{vol}_g$ is the Riemannian volume element (in local coordinates $x^\mu$,
$d\hbox{vol}_g = |g|dx$) and $|\partial u|_g^2:= g^{\mu\nu}\partial_\mu u\partial_\nu u$.
Similarly for any space-like hypersurface $\sigma$ we let $d\mu_g$ denote the positive 
Riemannian measure on $\sigma$ and $N$ be the future oriented unit normal vector to $\sigma$.
We then define
\[
 L^2(\sigma):= \{v: \sigma \longrightarrow \R \hbox{ measurable s.t. } 
 \|v\|_{L^2(\sigma)}^2:= \int_\sigma v^2d\mu_g <+\infty\},
\]
and, using a $g$-orthonormal frame $(e_0,\cdots,e_n)$ s.t. $e_0=N$,
\[
 H^1_m(\sigma):= \{v: \sigma \longrightarrow \R \hbox{ measurable s.t. }
\int_\sigma \left(\Sigma_{i=0}^n\langle e_i,\nabla u\rangle_g^2  +m^2u^2\right)d\mu_g <+\infty\}.
\]
A hypersurface $\Sigma$ is called
\emph{Cauchy} if any maximal smooth causal curve in $\mathcal{M}$ intersects $\Sigma$ at exactly one point
(a smooth causal curve is a a curve s.t. any vector which is tangent to it is time-like).
If $\Sigma_1$ and $\Sigma_2$ are two space-like hypersurfaces, we write
\emph{$\Sigma_1\prec \Sigma_2$ if $\Sigma_1$ is in the past of $\Sigma_2$ and $\Sigma_1\cap \Sigma_2 = \emptyset$}.
If $u$ is a real valued map defined on a neighbourhood of $\Sigma$, we denote by
$[u]_\Sigma =(u|_\Sigma,\langle N,\nabla u\rangle_g|_\Sigma)$ the Cauchy data of $u$
along $\Sigma$.
Our aim is to prove the existence and uniqueness of weak solutions to (\ref{KG-on-curved})
with Cauchy conditions in $H^1_m(\Sigma)\times L^2(\Sigma)$ for some space-like
Cauchy hypersurface $\Sigma$.

\subsection{Existence of solutions to the linear problem}
We first need generalizations of Proposition \ref{proposition-continuity-of-Phi} to this context. 
Such results were proved by Y. Choquet-Bruhat, D. Chistodoulou and M. Francaviglia \cite{cbcf}.
Here we present a more general version of their result in the case $s = r =1$ by using the same techniques
(see also the beautiful book by S. Alinhac \cite{alinhac}).
We will make the following further hypotheses on $(\mathcal{M},g)$:
there exists a smooth `temporal function' $\tau:\mathcal{M}\longrightarrow \R$ and
a smooth `radial function' $\rho:\mathcal{M}\longrightarrow [0,+\infty)$ and constants $A_1,A_2,A_3>0$ s.t.
\begin{equation}\label{hypo-tau-rho-1}
 |\nabla \tau|^2_g >0\quad \hbox{everywhere};
\end{equation}
\begin{equation}\label{hypo-tau-rho-2}
 \forall t\in \R,\quad
 \Sigma_t:= \tau^{-1}(t)\quad \hbox{is a space-like Cauchy hypersurface};
\end{equation}
\begin{equation}\label{hypo-tau-rho-3}
\forall r>0,\forall t_1,t_2\in \R\hbox{ s.t. }t_1<t_2,\quad
\{x\in \mathcal{M}; \tau(x)\in [t_1,t_2],\rho(x)\leq r\}
\hbox{ is compact}.
\end{equation}
Moreover there exists some $R_0>0$, s.t., $\forall x\in \mathcal{M}$,
\begin{equation}\label{hypo-tau-rho-4}
\rho(x)\geq R_0\Longrightarrow 
-A_3\leq |\nabla \rho(x)|^2_g <0
\quad\hbox{and}\quad |\langle \nabla \rho,\nabla \tau(x)\rangle_g| \leq A_2/\rho(x);
\end{equation}
\begin{equation}\label{hypo-tau-rho-5}
\rho(x)\geq R_0\Longrightarrow A_1/\rho(x)^2\leq |\nabla \tau(x)|^2_g.
\end{equation}
Lastly define $\ell:= |\nabla \tau|_g^{-1}$ (the \emph{lapse function}) and
$T:= \ell \nabla \tau = \nabla \tau /|\nabla \tau|_g$. We assume
that there exists a continuous function $B:\R \longrightarrow [0,+\infty)$ s.t.
\begin{equation}\label{further-hypo}
(n+1)\left(|\nabla T|_g^2 - 2|T^\mu\nabla_\mu T|_g^2\right)^{1/2}
\leq (B\circ \tau)|\nabla \tau|_g\quad \hbox{on }\mathcal{M}.
\end{equation}
Conditions (\ref{hypo-tau-rho-1}) and (\ref{hypo-tau-rho-2}) are equivalent to the assumption that
$(\mathcal{M},g)$ is \emph{globally hyperbolic} (see \cite{bgp}). 
Conditions (\ref{hypo-tau-rho-5}) (together with (\ref{hypo-tau-rho-4})) means that the
lapse function grows at most linearly in $\rho$ at spatial infinity. Condition (\ref{further-hypo}) is
an assumption on the curvature of the integral curves of the vector field $T$.

Given a smooth function $u:\mathcal{M}\longrightarrow \R$ we define its \emph{stress-energy tensor} $S(u)$
(associated with the action functional $\mathcal{A}$, see \cite{helein-cup}),
defined by
\[
 S^\mu_\nu(u):= g^{\mu\lambda}\partial_\lambda u\partial_\nu u -
 \frac{1}{2}\left(|\partial u|_g^2 - m^2u^2\right)\delta^\mu_\nu.
\]
We say that \emph{$u$ has a compact spatial support} if it vanishes on
$\{x\in \mathcal{M}; \rho(x)\geq h(t)\}$ for some continuous function $h:\R\longrightarrow [0,+\infty)$.
If so and if $\sigma$ is a space-like hypersurface (possibly with boundary), we define the energy
\[
 E_u(\sigma):= \int_\sigma\langle S(u) N,N\rangle_g d\mu_g.
\]
Note that $E_u(\sigma)$ is always nonnegative. In particular
if, on $\sigma$, we use a $g$-orthonormal frame $(e_0,\cdots,e_n)$ s.t. $e_0=N$,
then $\langle S(u) N,N\rangle_g = S^0_0(u)= \frac{1}{2}\left(
\Sigma_{\mu=0}^n\langle e_\mu,\nabla u\rangle_g^2  +m^2u^2\right)$.

For any interval $I \subset \R$, we define $\|u\|_{I,\tau}:= \sup_{t\in I} E_u(\Sigma_t)^{1/2}$ and
\[
 \mathcal{F}^1_I(\Sigma_\tau):= \left\{\begin{array}{l}
                                            \hbox{the closure of the set of smooth functions on }\tau^{-1}(I)\hbox{ with}\\
                                            \hbox{compact spatial support in the topology induced by }\|\cdot\|_{I,\tau}.
                                           \end{array}
 \right.
\]
We also define
$\mathcal{F}^1_{loc}(\Sigma_\tau):= \{u:\mathcal{M}\longrightarrow \R; \forall I\subset \R\hbox{ s.t. }I\hbox{ is bounded},
\|u\|_{I,\tau} <+\infty\}$.
For any positive function $\beta:\R\longrightarrow (0,+\infty)$ (actually we will use
$\beta(t) = \hbox{exp}\frac{1}{2}\int_0^t|B(s)|ds$), we set
$\|u\|_{\beta,\tau}:= \sup_{t\in \R} \beta(t)^{-1}E_u(\Sigma_t)^{1/2}$
and $\mathcal{F}^1_{\beta}(\Sigma_\tau):= \{u\in \mathcal{F}^1_{loc}(\Sigma_\tau);\|u\|_{\beta,\tau} <+\infty\}$. 
Lastly we set
\[
 L^1_{loc}(\R,L^2_\ell(\Sigma_\tau)):= \{f:\mathcal{M}\longrightarrow \R\hbox{ measurable s.t. }
 [t\longmapsto \|\ell f|_{\Sigma_t}\|_{L^2(\Sigma_t)}]\in L^1_{loc}(\R)\}
\]
and, for any interval $I\subset \R$, we note $\|f\|_{L^1(I,L^2_\ell(\Sigma_\tau))}:= 
\int_I \|\ell f\|_{L^2(\Sigma_t)}dt$.

The existence result in \cite{cbcf} concerned weak solutions to (\ref{KG-on-curved})
with a Cauchy data on a hypersurface $\Sigma_t$. The following result extends this with
the notable difference that we allow more general Cauchy hypersurfaces.  Fixing $\tau$ (and hence the foliation
$(\Sigma_t)_t$) we say that \emph{a space-like hypersurface $\widehat{\Sigma}$ is admissible if it is
a Cauchy hypersurface and if:
(i) $\exists t_1,t_2\in \R$ s.t. $t_1<t_2$ and $\Sigma_{t_1} \prec \widehat{\Sigma} \prec \Sigma_{t_2}$;
(ii) if $\widehat{N}$ denotes the future oriented normal to $\widehat{\Sigma}$,
$C(\widehat{\Sigma}):= \sup_{\widehat{\Sigma}}\langle \widehat{N},T\rangle_g<+\infty$}.
\begin{theo}\label{existence-on-curved}
Assume that $(\mathcal{M},g)$ satisfies Hypotheses (\ref{hypo-tau-rho-1})--(\ref{further-hypo}).
Let $\widehat{\Sigma}\subset \mathcal{M}$ be a space-like admissible hypersurface
and let $t_1,t_2\in \R$ s.t. $\Sigma_{t_1} \prec \widehat{\Sigma} \prec \Sigma_{t_2}$.
Then for any $(u_0,u_1)\in H^1_m(\widehat{\Sigma})\times  L^2(\widehat{\Sigma})$
and $f\in L^1_{loc}(\R,L^2_\ell(\Sigma_\tau))$, there
exists an unique weak solution $u\in \mathcal{F}^1_{loc}(\Sigma_\tau)$ to (\ref{KG-on-curved}) 
s.t. $[u]_{\widehat{\Sigma}} = (u_0, u_1)$.
Moreover for $t = t_1$ or $t_2$,
\begin{equation}\label{main-estimate-on-curved}
 E_u(\Sigma_{t})^{1/2} \leq  \sqrt{2}\beta(t_1,t_2)\left(C(\widehat{\Sigma})
 E_u(\widehat{\Sigma})\right)^{1/2}
 + \sqrt{2}\beta(t_1,t_2)^2 \|f\|_{L^1([t_1,t_2],L^2_\ell(\Sigma_\tau))},
\end{equation}
where $\beta(t_1,t_2):= e^{\frac{1}{2}|\int_{t_1}^{t_2} B(s)ds|}$ and
$C(\widehat{\Sigma}):= \sup_{\widehat{\Sigma}}\langle \widehat{N},T\rangle_g$.
\end{theo}
\emph{Proof} --- The main point is to obtain the a priori estimate (\ref{main-estimate-on-curved})
for any solution $u$ to (\ref{KG-on-curved}). 
Without loss of generality we will content ourself to prove that
\begin{equation}\label{simpli-estimate-on-curved}
 E_u(\Sigma_{t_2})^{1/2} \leq \sqrt{2}
 e^{\frac{1}{2}\int_{t_1}^{t_2} B(s)ds}\left( C(\widehat{\Sigma})E_u(\widehat{\Sigma})\right)^{1/2} +
 \sqrt{2}e^{\int_{t_1}^{t_2} B(s)ds}\int_{t_1}^{t_2} \left\|\ell f|_{\Sigma_s}\right\|_{L^2}ds.
\end{equation}

\noindent
\emph{Step 1: Use of conservation law} --- 
Consider a compact domain $D\subset \tau^{-1}([t_1,t_2])$,
the boundary $\partial D$ of which is composed of three smooth components
\[
 \partial D = (\Sigma_{t_2}\cap \overline{D}) + \Lambda - (\Sigma_{t_1}\cap \overline{D}),
\]
where the signs give the orientation.
We assume that $\Lambda$ is space-like and that the normal vector
$N$ to it is future-pointing, hence $(\Sigma_{t_2}\cap \overline{D}) \cup \Lambda$ forms the top of
$D$, whereas $\Sigma_{t_1}\cap \overline{D}$ is the bottom (see the end of the proof
for the construction of $D$).

Let $\theta:\mathcal{M}\longrightarrow \R$
be the function which coincides with $\tau$ on $\widehat{\Sigma}$ and which is
invariant by the flow of $T$ and, for $t\in [t_1,t_2]$, consider the domain
\[
 D_t:= \{x\in D; \theta(x)<\tau(x)<t\}
\]
(points in $D_t$ are points of $D$ which are in the future of $\widehat{\Sigma}$ and in the past
of $\Sigma_t$, see the figure).
\begin{figure}[h]
\begin{center}
\input{sigma.pstex_t}
\end{center}
\end{figure}
Note that $\partial D_t = (\Sigma_t\cap \overline{D_t}) + (\Lambda \cap \overline{D_t})
- \widehat{\Sigma}_{<t}$, where
$\Sigma_t\cap \overline{D_t} = \{ x\in \overline{D}; \theta(x)<\tau(x) = t\}$
and $\widehat{\Sigma}_{<t} = \{ x\in \overline{D}; \theta(x)=\tau(x)<t\}$.

Let us apply Stokes theorem to 
$S(u)T$ on $D_t$. We get (writing $S = S(u)$ for shortness):
\begin{equation}\label{Stokes-de-depart}
 \int_{D_t} \nabla_\mu (S^\mu_\nu T^\nu) d\hbox{vol}_g
 = \int_{\partial D_t} \langle N,ST\rangle_g d\mu_g.
\end{equation}
Since $u$ is a solution of (\ref{KG-on-curved}), the stress-energy tensor satisfies the relation
$\nabla_\mu S^\mu_\nu = f\partial_\nu u$, see e.g. \cite{helein-cup}.
Hence the l.h.s. of (\ref{Stokes-de-depart}) reads
\begin{equation}\label{int-nabla-Y}
\int_{D_t} \nabla_\mu (S^\mu_\nu T^\nu)d\hbox{vol}_g 
 =  \int_{D_t} f (\partial_\nu u) T^\nu d\hbox{vol}_g
 + \int_{D_t} S^\mu_\nu\nabla_\mu  T^\nu d\hbox{vol}_g.
\end{equation}
(i) \emph{Estimation of the first term in the r.h.s of (\ref{int-nabla-Y})} --- Using the coarea formula, we get
\[
 \left| \int_{D_t} f T^\nu \partial_\nu u\, d\hbox{vol}_g \right|
 \leq \int_{D_t} |f| |T^\nu \partial_\nu u| d\hbox{vol}_g \\
 = \int_{t_1}^tds\int_{\Sigma_s\cap \overline{D_t}}|T^\nu \partial_\nu u| |f| \frac{d\mu_g}{|\nabla \tau|_g}.
\]
Since $T$ coincides with the normal vector $N$ to $\Sigma_s$, we have:
$|T^\nu \partial_\nu u| = |\langle N,\nabla u\rangle_g|\leq \sqrt{2\langle S N,N\rangle_g}$.
Hence by Cauchy--Schwarz
\[
\begin{array}{ccl}
 \displaystyle \left| \int_{D_t} f T^\nu \partial_\nu u d\hbox{vol}_g \right| & \leq &
 \displaystyle \int_{t_1}^tds \left(\int_{\Sigma_s\cap \overline{D_t}}2\langle S N,N\rangle_gd\mu_g\right)^{1/2}
 \left(\int_{\Sigma_s\cap \overline{D_t}}  (\ell f)^2d\mu_g\right)^{1/2}\\
 & = & \displaystyle \int_{t_1}^tds \sqrt{2E_u(\Sigma_s\cap \overline{D_t})}
 \left\|\ell f|_{\Sigma_s\cap \overline{D_t}}\right\|_{L^2}.
\end{array}
\]
(ii) \emph{Estimation of the second term in the r.h.s of (\ref{int-nabla-Y})} --- The 
Cauchy--Schwarz inequality gives us:
\[
 |S^\mu_\nu\nabla_\mu T^\nu|^2  \leq \left(\sum_{\mu,\nu=0}^n(S^\mu_\nu)^2\right)\left(\sum_{\mu,\nu=0}^n(\nabla_\mu T^\nu)^2\right),
\]
however a difficulty is that the r.h.s. of this inequality depends on the choice of the frame $(e_0,\cdots,e_n)$ used
in the decomposition of the tensors $S$ and $\nabla T$. We choose a $g$-orthonormal frame $(e_0,\cdots,e_n)$ s.t. $e_0=T$.
Observe then that $|S^\mu_\nu|\leq S^0_0$, $\forall \mu,\nu$ and thus\footnote{This inequality is true for any vector
valued field $u$.  Actually using the fact that
$u$ is a scalar field one can get the improved inequality $\sum_{\mu,\nu=0}^n(S^\mu_\nu)^2 \leq (n+3)(S^0_0)^2$.}
\[
 \sum_{\mu,\nu=0}^n(S^\mu_\nu)^2 \leq (n+1)^2(S^0_0)^2 = (n+1)^2\langle T,ST\rangle_g^2.
\]
Next let us introduce the tensor $h_{\mu\nu}:= 2T_\mu T_\nu - g_{\mu\nu}$.
We note that in the previously chosen $g$-orthonormal frame we have $h_{\mu\nu} = \delta_{\mu\nu}$
and hence that
\[
 \sum_{\mu,\nu=0}^n(\nabla_\mu T^\nu)^2 = \nabla_\mu T^\nu \nabla_\lambda T^\sigma h_{\nu\sigma}h^{\mu\lambda}
 =: |\nabla T|_h^2.
\]
We thus deduce that
\[
 |S^\mu_\nu\nabla_\mu T^\nu| \leq (n+1)\langle T,ST\rangle_g |\nabla T|_h,
\]
where the r.h.s. is now frame independent. Lastly a computation
(using $|T|_g^2=1$, which implies $g_{\lambda\nu}T^\nu\nabla_\mu T^\lambda = 0$) shows that
$|\nabla T|_h^2 = |\nabla T|^2_g -2|T^\mu\nabla_\mu T|_g^2$.
Hence $|S^\mu_\nu\nabla_\mu T^\nu| \leq (n+1)(|\nabla T|^2_g -2|T^\mu\nabla_\mu T|_g^2)^{1/2}
 S^0_0$ and using the fact that $T$ coincides with the future pointing normal vector to $\Sigma_s$,
(\ref{further-hypo}) and the coarea formula,
 \[
 \begin{array}{ccl}
  \displaystyle 
  \left|\int_{D_t} S^\mu_\nu\nabla_\mu  T^\nu d\hbox{vol}_g\right| & \leq &
  \displaystyle \int_{D_t} B(\tau)|\nabla \tau|_gS^0_0  d\hbox{vol}_g\\
  & = & \displaystyle  \int_{t_1}^t B(s)ds\int_{\Sigma_s\cap \overline{D_t}}\langle ST,T\rangle_g d\mu_g
  = \int_{t_1}^t B(s) E_u(\Sigma_s\cap \overline{D_t})ds .
 \end{array}
 \]
Summarizing with the previous step we deduce from (\ref{int-nabla-Y})
\[
 \left|
\int_{D_t} \nabla_\mu (S^\mu_\nu T^\nu)d\hbox{vol}_g\right| \leq
\int_{t_1}^t \left\|\ell f|_{\Sigma_s}\right\|_{L^2}\sqrt{2E_u(\Sigma_s\cap \overline{D_t})}ds
+ \int_{t_1}^t B(s) E_u(\Sigma_s\cap \overline{D_t})ds,
\]
which, in view of (\ref{Stokes-de-depart}) gives:
\begin{equation}\label{estimate-on-left-Stokes}
 \int_{\partial D_t} \langle N,ST\rangle_g d\mu_g \leq
\int_{t_1}^t \left\|\ell f|_{\Sigma_s}\right\|_{L^2}\sqrt{2E_u(\Sigma_s\cap \overline{D_t})}ds
+ \int_{t_1}^t B(s) E_u(\Sigma_s\cap \overline{D_t})ds.
\end{equation}
(iii) \emph{Lower estimation of the l.h.s. of (\ref{estimate-on-left-Stokes})} ---
Using the fact that $T=N$ on $\Sigma_t$ and denoting by $\widehat{N}$ the future pointing
normal to $\widehat{\Sigma}$, we decompose
\[ \int_{\partial D_t} \langle N,ST\rangle_g d\mu_g
  = E_u(\Sigma_t\cap \overline{D_t})
 + \int_{D_t\cap \Lambda}\langle ST,N\rangle_g d\mu_g
 - \int_{\widehat{\Sigma}_{<t}}\langle ST,\widehat{N}\rangle_gd\mu_g.
\]
However $\langle ST,N\rangle_g \geq 0$ on $D_t\cap \Lambda$. This follows from
$|N|_g^2 = 1$, $N^0>0$ and from the following identity, valid in a $g$-orthonormal frame $(e_0,\cdots,e_n)$ s.t. $e_0=T$:
\begin{equation}\label{identity-on-stress-energy}
2N^0\langle ST,N\rangle_g = |N|_g^2(u^0)^2 + \sum_{i=1}^n(N^0u^i-N^iu^0)^2 + (N^0)^2m^2u^2,
\end{equation}
where $N^\mu = \langle e_\mu,N\rangle_g$ and $u^\mu:= \langle e_\mu,\nabla u\rangle_g$. Hence
\begin{equation}\label{conclusion-de-Stokes}
  E_u(\Sigma_t\cap \overline{D_t})
 - \int_{\widehat{\Sigma}_{<t}}\langle ST,\widehat{N}\rangle_gd\mu_g
 \leq \int_{\partial D_t} \langle N,ST\rangle_g d\mu_g.
\end{equation}
(iv) \emph{Conclusion} --- For any $s\in [t_1,t]$
let $\widehat{\Sigma}_{\geq s}:= \{x\in \widehat{\Sigma}; \tau(x)\geq s\}$
and set
\[
 e(s):= E_u(\Sigma_s\cap \overline{D_t}) + \int_{\widehat{\Sigma}_{\geq s}}\langle ST,\widehat{N}\rangle_gd\mu_g.
\]
We will prove that the l.h.s. of (\ref{conclusion-de-Stokes}) is equal to $e(t)-e(t_1)$.
Observe that, because of $\Sigma_{t_1}\cap \overline{D_t} = \emptyset$ and
$\widehat{\Sigma}_{\geq t_1} = \widehat{\Sigma}\cap D$,
$e(t_1) = \int_{\widehat{\Sigma}\cap D}\langle ST,\widehat{N}\rangle_g d\mu_g$.
But, since $\widehat{\Sigma}_{\geq t}\cap \widehat{\Sigma}_{<t} = \emptyset$ and
$\widehat{\Sigma}_{\geq t}\cup \widehat{\Sigma}_{<t} = \widehat{\Sigma}\cap D$, the latter decomposes as:
\[
  e(t_1) = \int_{\widehat{\Sigma}\cap D}\langle ST,\widehat{N}\rangle_gd\mu_g
  = \int_{\widehat{\Sigma}_{<t}}\langle ST,\widehat{N}\rangle_gd\mu_g
 + \int_{\widehat{\Sigma}_{\geq t}}\langle ST,\widehat{N}\rangle_gd\mu_g.
\]
Hence $e(t) - e(t_1) = E_u(\Sigma_t\cap\overline{D_t}) - \int_{\widehat{\Sigma}_{<t}}\langle ST,\widehat{N}\rangle_gd\mu_g$,
so that (\ref{conclusion-de-Stokes}) reads
$e(t) - e(t_1) \leq \int_{\partial D_t} \langle N,ST\rangle_g d\mu_g$.
By using (\ref{estimate-on-left-Stokes}) and the fact that
$E_u(\Sigma_s\cap \overline{D_t})\leq e(s)$ we deduce 
(setting $F(s):= \sqrt{2}\left\|\ell f|_{\Sigma_s}\right\|_{L^2}$) that:
\begin{equation}\label{pre-gronwall}
 e(t) - e(t_1)\leq \int_{t_1}^tds F(s)\sqrt{e(s)}
+ \int_{t_1}^t B(s) e(s)ds.
\end{equation}
\emph{Step 3: Using Gronwall lemma} --- Set $K:= e(t_1) + \int_{t_1}^tds F(s)\sqrt{e(s)}$.
Then (\ref{pre-gronwall}) (by replacing $t$ by $t'$) implies easily
$e(t')\leq K + \int_{t_1}^{t'} B(s)e(s)ds,\quad \forall t'\in [t_1,t]$.
Using Gronwall Lemma we deduce that
\[
e(t')\leq K e^{\int_{t_1}^{t'} B(s)ds}, \quad \forall t'\in [t_1,t].
\]
Replacing $K$ by its value, setting $\psi(t):= \sup_{t_1\leq s\leq t} \sqrt{e(s)}$ and taking the supremum
over $t'\in [t_1,t]$, we obtain
\[
 \psi(t)^2\leq \left(\psi(t_1)^2 + \int_{t_1}^tds F(s)\psi(s) \right)e^{\int_{t_1}^t B(s)ds}
 \leq \left(\psi(t_1)^2 + \psi(t)\int_{t_1}^t F(s)ds \right)e^{\int_{t_1}^t B(s)ds},
\]
which implies $\psi(t) \leq \psi(t_1)e^{\frac{1}{2}\int_0^t B(s)ds} + \left(\int_{t_1}^t F(s)ds\right)e^{\int_{t_1}^t B(s)ds}$.
Applying this for $t=t_2$ and using $e(t_2) = E_u(\Sigma_{t_2}\cap \overline{D})$, we get
\begin{equation}\label{intermediate-energy-estimate}
E_u(\Sigma_{t_2}\cap \overline{D})^{1/2} \leq e^{\frac{1}{2}\int_{t_1}^{t_2} B(s)ds}e(t_1)^{1/2}
+ e^{\int_{t_1}^{t_2} B(s)ds}\int_{t_1}^{t_2} F(s)ds.
\end{equation}
\emph{Step 4: Controlling $e(t_1)$ by $E_u(\widehat{\Sigma})$} ---
Using an identity similar to (\ref{identity-on-stress-energy}) (where $T$
is replaced by $\widehat{N}$, $N$ is replaced by $T$ and we use a $g$-orthonormal
frame $(\widehat{e}_0,\cdots,\widehat{e}_n)$ s.t. $\widehat{e}_0 = \widehat{N}$)
we prove that
\[
 \langle ST,\widehat{N}\rangle_g \leq 2\langle T,\widehat{N}\rangle_g
 \left(\frac{1}{2}\sum_{\mu=0}^n\langle \widehat{e}_\mu,\nabla u\rangle_g^2 + \frac{1}{2}m^2u^2\right)
 = 2\langle T,\widehat{N}\rangle_g\langle S\widehat{N},\widehat{N}\rangle_g.
\]
This hence implies that
\[
 e(t_1) = \int_{\widehat{\Sigma}\cap D}\langle ST,\widehat{N}\rangle_g d\mu_g
 \leq 2 \sup_{\widehat{\Sigma}}\langle T,\widehat{N}\rangle_gE_u(\widehat{\Sigma}\cap \overline{D})
 \leq 2 C(\widehat{\Sigma})E_u(\widehat{\Sigma}).
\]
Thus we deduce from (\ref{intermediate-energy-estimate})
\begin{equation}\label{intermediate-energy-estimate2}
E_u(\Sigma_{t_2}\cap \overline{D})^{1/2} \leq \sqrt{2}e^{\frac{1}{2}\int_{t_1}^{t_2} B(s)ds}
\left(C(\widehat{\Sigma})E_u(\widehat{\Sigma})\right)^{1/2}
+ e^{\int_{t_1}^{t_2} B(s)ds}\int_{t_1}^{t_2} F(s)ds.
\end{equation}
\emph{Step 5: Global estimate} ---
Now, for any $R>0$, set $K_R:= \{x\in \mathcal{M}; \rho(x)\leq R\}$.
In order to obtain (\ref{simpli-estimate-on-curved}) it suffices to prove that there exists
some $R_0>0$ s.t., for any $R>R_0$, there exists a domain $D$ satisfying the previous properties and s.t.
$\Sigma_{t_2}\cap K_R\subset \Sigma_{t_2}\cap \overline{D}$. Indeed if so we deduce from (\ref{intermediate-energy-estimate2})
\[
  E_u(\Sigma_{t_2}\cap K_R)^{1/2} \leq E_u(\Sigma_{t_2}\cap \overline{D})^{1/2}  \leq 
 \sqrt{2} e^{\frac{1}{2}\int_{t_1}^{t_2} B(s)ds}\left(C(\widehat{\Sigma})E_u(\widehat{\Sigma})\right)^{1/2}
 + e^{\int_{t_1}^{t_2} B(s)ds}\int_{t_1}^{t_2} F(s)ds
\]
Since this inequality holds for any $R>0$, it thus implies (\ref{simpli-estimate-on-curved}).

\noindent
\emph{Step 6: Construction of $D$} --- Here we need Hypotheses (\ref{hypo-tau-rho-3})
to (\ref{hypo-tau-rho-5}). Set $t:= t_2-t_1$ and assume that $t>0$.
For any $R>R_0$ we will construct a smooth function
$\tilde{\tau}:\mathcal{M}\longrightarrow \R$ and find some $\overline{R}> R$ s.t.
\begin{enumerate}
 \item $\forall x\in K_R$, $\tilde{\tau}(x) = \tau(x) - t = \tau(x) -t_2+t_1$;
 \item $\forall x\not\in K_{\overline{R}}$, $\tilde{\tau}(x) = \tau(x)$;
 \item $|\nabla \tilde{\tau}|_g^2 >0$ everywhere, in particular the level sets of $\tilde{ \tau}$ are space-like
hypersurfaces.
\end{enumerate}
If so $D: = \{x\in \mathcal{M}; \tau(x)>t_1,\tilde{\tau}(x)<t_1\}$ satisfies all the previously required properties.
To construct $\tilde{ \tau}$, we set $\tilde{ \tau} = \tau - t\chi\circ \rho$, where $\chi\in \mathcal{C}^0([0,+\infty),[0,1])$
is piecewise $\mathcal{C}^\infty$ and has to be suitably chosen. Conditions (i) and (ii) translate
respectively as:
(i)' $\forall r\leq R$, $\chi(r) = 1$;
(ii)' $\forall r\geq \overline{R}$, $\chi(r) = 0$.
A simple computation using (\ref{hypo-tau-rho-4})
and (\ref{hypo-tau-rho-5}) shows that Condition (iii) is satisfied if
\[
A_3t^2r^2(\chi'(r))^2 + 2A_2tr|\chi'(r)| < A_1.
\]
This condition is fulfilled if we choose $\alpha>0$ s.t. $A_3\alpha^2+2A_2\alpha < A_1$, $\overline{R} = Re^{t/\alpha}$
and set $\chi(r) = 1 - \frac{\alpha}{t} \log\frac{r}{R}$, $\forall r\in [R,\overline{R}]$
(all that works because $\int_R^\infty \frac{dr}{r} = + \infty$).\\

\noindent
\emph{Step 7: Conclusion} --- Thanks to the works of J. Hadamard, M. Riesz and the results by J. Leray \cite{leray},
one can construct fundamental
solutions for the operator $L$ and solve the Cauchy problem for smooth Cauchy data (see
\cite{friedlander,bgp}). By using the density of smooth compactly supported functions
in $L^2(\Sigma_{t_1})$ and $H^1_m(\Sigma_{t_1})$ and (\ref{main-estimate-on-curved}), we deduce the existence.
The uniqueness is a straightforward consequence of (\ref{main-estimate-on-curved}).
\hfill $\square$\\

\noindent
Note that similar results exist for higher (integer) order Sobolev spaces and for
Cauchy data on a hypersurface which belongs to the family $(\Sigma_t)_t$, see
\cite{cbcf} and also \cite{alinhac}. Theorem \ref{existence-on-curved} has the following
consequence which is a substitute for Proposition \ref{proposition-continuity-of-Phi}.
Set $\mathcal{E}^1_{0,\beta}(\Sigma_\tau):= \{\varphi\in \mathcal{F}^1_\beta(\Sigma_\tau));
\square_g\varphi + m^2\varphi = 0\}$. 
\begin{coro}\label{coro-continuity-of-Phi-curved}
Let $(\mathcal{M},g)$ be a $n$-dimensional Lorentzian manifold. Assume that there exist
functions $\tau,\rho\in \mathcal{C}^\infty(\mathcal{M},\R)$ which satisfy
(\ref{hypo-tau-rho-1})--(\ref{further-hypo}). Set $\beta(t) = e^{\frac{1}{2}|\int_0^t B(s)ds|}$.
Then for any admissible hypersurface $\sigma$, there exists a continuous linear map
\[
\begin{array}{cccl}
 \Phi_\sigma : & H^1_m(\sigma)\times L^2(\sigma) & \longrightarrow & \mathcal{E}^1_{0,\beta}(\Sigma_\tau)\\
 & (\psi,\chi) & \longmapsto & \Phi_\sigma(\psi,\chi),
\end{array}
\]
where $\Phi_\sigma(\psi,\chi)$ is equal to the unique solution $\varphi$ to $L_g\varphi = \square_g\varphi + m^2\varphi = 0$
with the Cauchy data $[\varphi]_\sigma = (\psi,\chi)$.
\end{coro}
Thanks to this result we can define for any admissible hypersurface $\sigma$
the continuous map
\[
 \Theta_\sigma: \mathcal{F}^1_\beta(\Sigma_\tau)\longrightarrow \mathcal{E}^1_{0,\beta}(\Sigma_\tau)
\]
defined by $\Theta_\sigma(u) := \Phi_{\sigma}([u]_{\sigma})$.
The following result will also be useful.
\begin{lemm}\label{lemma-de-Sigma0-vers-tildeSigma} 
 Let $f\in L^1_{loc}(\R,L^2_\ell(\Sigma_\tau))$ and
 $u\in \mathcal{F}^1_{loc}(\Sigma_\tau)$ be a solution of $\square_gu+m^2u = f$. Let
 $\widehat{\Sigma}$ be an admissible hypersurface s.t. $\Sigma_{t_1} \prec \widehat{\Sigma} \prec\Sigma_{t_2}$.
 Then
 \begin{equation}\label{estimate-de-Sigma0-vers-tildeSigma}
  E_u(\widehat{\Sigma}) \leq 2C(\widehat{\Sigma}) \left[
  \left(1+\|B\|_{L^1([t_1,t_2])}\right) \|u\|_{[t_1,t_2],\tau}
  + \sqrt{2}\|f\|_{L^1([t_1,t_2],L_\ell^2)}\|u\|_{[t_1,t_2],\tau}^{1/2}\right].
 \end{equation}
\end{lemm}
\emph{Sketch of the proof} --- The proof is based on the same techniques as in the proof
of Theorem \ref{existence-on-curved}: one starts from the identity
$\int_{\Delta} \nabla_\mu(S^\mu_\nu T^\nu)d\hbox{vol}_g =
\int_{\partial \Delta} \langle ST,N\rangle_g d\mu_g$,
with the same vector field $T$. The difference is the domain of integration which is now
$\Delta:= \{x\in D; t_1<\tau(x)<\theta(x)\}$. Also the reasoning is simpler, for 
we already know that $\|u\|_{[t_1,t_2],\tau}$ is bounded and hence 
we do not need to use Gronwall lemma. This leads to
\[
 \int_{\widehat{\Sigma}} \langle ST,\widehat{N}\rangle_g d\mu_g \leq
 \left(1+\|B\|_{L^1([t_1,t_2])}\right) \|u\|_{[t_1,t_2],\tau}
  + \sqrt{2}\|f\|_{L^1([t_1,t_2],L_\ell^2)}\|u\|_{[t_1,t_2],\tau}^{1/2}.
\]
Estimate (\ref{estimate-de-Sigma0-vers-tildeSigma}) follows then from the inequality
$\langle S\widehat{N},\widehat{N}\rangle_g \leq 2\langle T,\widehat{N}\rangle_g \langle ST,\widehat{N}\rangle_g$,
which implies $E_u(\widehat{\Sigma}) \leq 2\sup_{\widehat{\Sigma}} \langle T,\widehat{N}\rangle_g
\int_{\widehat{\Sigma}} \langle ST,\widehat{N}\rangle_g d\mu_g$.\hfill $\square$

\subsection{A generalization of Duhamel's formula}\label{subsec-gen-Duhamel}
Our aim is here to prove a `curved' version of Duhamel's formula.
Beside the foliation of $\mathcal{M}$ by the level sets $\Sigma_t:= \tau^{-1}(t)$,
we also consider a family $(\sigma_s)_{s\in \R}$ of admissible Cauchy
space-like hypersurfaces, \emph{which may not form a foliation of $\mathcal{M}$
in general}. We assume that there exists an $n$-dimensional manifold $\underline{\sigma}$
(the model for each $\sigma_s$) and a map $F\in \mathcal{C}^\infty(\R\times \underline{\sigma},\mathcal{M})$
s.t. for any $s\in \R$, $F_s:= F(s,\cdot)$ is an embedding of $\underline{\sigma}$, the image
of which is $\sigma_s$.
On each $\sigma_s$ we define the function $\lambda_s\in  \mathcal{C}^\infty(\sigma_s,\R)$ by
$\lambda_s\circ F_s:= \langle \frac{\partial F}{\partial s}(s,\cdot),N_s\circ F_s\rangle_g$, where
$N_s$ is the future pointing normal vector to $\sigma_s$.
We call $(\sigma_s)_{s\in \R}$ a \emph{smooth family of admissible Cauchy hypersurfaces}.

For any $s\in \R$, we denote by $\{x\succ \sigma_s\}$ (resp. $\{x\prec \sigma_s\}$)
the subset of $\mathcal{M}\setminus \sigma_s$
which are in the future (resp. the past) of $\sigma_s$, similarly
$\{x\succcurlyeq \sigma_s\} := \sigma_s\cup \{x\succ \sigma_s\}$
($\{x\preccurlyeq \sigma_s\} := \sigma_s\cup \{x\prec \sigma_s\}$).
We let $Y_{\sigma_s}\in L^\infty(\mathcal{M})$ be s.t. $Y_{\sigma_s}=1$ on $\{x\succcurlyeq \sigma_s\}$
and $Y_{\sigma_s}=0$ on $\{x\prec \sigma_s\}$. 

We let $f\in L^1_{loc}(\R,L^2_\ell(\Sigma_\tau))$ and we assume that,
for a.e. $s\in \R$, $\lambda_sf|_{\sigma_s}\in L^2(\sigma_s)$ and
$[s\longmapsto \|\lambda_sf|_{\sigma_s}\|_{L^2}]$ belongs to $L^1_{loc}(\R)$.
We then define:
\[
 \gamma_sf:\quad \hbox{the unique solution of}\quad
 \left\{ \begin{array}{cccl}
          \gamma_sf & = & 0 & \hbox{on }\{x\prec \sigma_s\}\\
          \gamma_sf  & = & \Phi_{\sigma_s}(0,\lambda_sf|_{\sigma_s}) & \hbox{on }\{x\succcurlyeq \sigma_s\},
         \end{array} \right.
\]
\[
 \Gamma_sf:\quad \hbox{the unique solution of}\quad
 \left\{ \begin{array}{cccl}
          \Gamma_sf & = & 0 & \hbox{on }\{x\prec \sigma_s\}\\
          L_g(\Gamma_sf)  & = & fY_{\sigma_s} & \hbox{on }\mathcal{M}.
         \end{array} \right.
\]
For any $y\in \mathcal{M}$ we let $G_y$ be the solution of $L_gG_y=0$ with the Cauchy data
$G_y|_{\sigma} = 0$ and $\langle N,\nabla G_y\rangle_g|_{\sigma} =\delta_y$, where $\sigma$ is
a Cauchy hypersurface which contains $y$. Then, still if $y\in \sigma$,
$Y_\sigma G_y$ is the retarded Green function for $L_g$ with 
source $\delta_y$ (see \cite{bgp}) for its existence).
Thus if $f$ is smooth, then we have the representation formulas
$(\gamma_sf)(x) = \int_{\sigma_s}f(y)(Y_{\sigma_s} G_y)(x)\lambda_s(y)d\mu_g(y)$
and
$(\Gamma_sf)(x) = \int_{\{y\succ \sigma_s\}}f(y)(Y_{\sigma_s} G_y)(x)d\hbox{vol}_g(y)$.
\begin{prop}\label{duhamel-on-curved}
 Let $(\mathcal{M},g)$ be a $n$-dimensional Lorentzian manifold. Assume that there exist a
 temporal function $\tau$ and a radial function $\rho$ which satisfy (\ref{hypo-tau-rho-1})--(\ref{further-hypo}).
 Let $(\sigma_s)_{s\in \R}$ be a $\underline{\sigma}$-family of
 admissible Cauchy hypersurfaces s.t.
 $\lim_{s\rightarrow +\infty}\left(\inf_{x\in \sigma_s}\tau(x)\right) = +\infty$.
 Let $f\in L^1_{loc}(\R,L^2_\ell(\Sigma_\tau))$
 s.t. $[s\longmapsto \|\lambda_sf|_{\sigma_s}\|_{L^2}]\in L^1_{loc}(\R)$.
 
 Then for any $u\in \mathcal{F}_{loc}^1(\Sigma_\tau)$ s.t. $L_gu= f$, we have, for any $s\in \R$,
\begin{equation}\label{duhamel-abstrait}
 u = \Theta_{\sigma_s}u + \Gamma_sf\quad \hbox{on }\{x\succ \sigma_s\}.
\end{equation}
Moreover 
\begin{equation}\label{B=intbeta}
 \Gamma_sf = \int_s^\infty (\gamma_{s_1}f)ds_1 .
\end{equation}
\end{prop}
\emph{Remark} --- The integral in the r.h.s. of (\ref{B=intbeta}) makes
sense as a distribution on $\mathcal{M}$ since, for any $\varphi\in \mathcal{C}^\infty_c(\mathcal{M})$,
we can set
$\langle \int_s^\infty (\gamma_{s_1}f)ds_1,\varphi\rangle = \int_s^\infty \langle \gamma_{s_1}f,\varphi\rangle ds_1
= \int_s^{\overline{s}} \langle \gamma_{s_1}f,\varphi\rangle ds_1$, where $\overline{s}$ is s.t.
$\hbox{supp}\varphi\subset \{ \sigma_s \prec x\prec \sigma_{\overline{s}}\}$
($\overline{s}$ exists because $\lim_{s\rightarrow +\infty}\left(\inf_{x\in \sigma_s}\tau(x)\right) = +\infty$).\\
\emph{Proof} --- The proof of (\ref{duhamel-abstrait}) is easy: since $\sigma_s$ is admissible,
there exists some $t\in \R$ s.t. $\Sigma_t \prec \sigma_s$ and thus $[\Gamma_sf]_{\Sigma_t}=0$.
Using arguments similar to the ones
used in the proofs of Theorem \ref{existence-on-curved} or Lemma \ref{lemma-de-Sigma0-vers-tildeSigma},
one can deduce that $E_{\Gamma_sf}(\sigma_s) = 0$, i.e. $[\Gamma_sf]_{\sigma_s}=0$. Hence the Cauchy data on $\sigma_s$
of both sides of (\ref{duhamel-abstrait}) coincide. Since these both sides are also solution
of the equation $L_g\varphi = f$ on $\{x\succ \sigma_s\}$, (\ref{duhamel-abstrait})
follows by uniqueness of the solution.

To prove (\ref{B=intbeta}), fix $s\in \R$ and set $v:= \int_s^\infty (\gamma_{s_1}f)ds_1$.
We take any $\varphi\in \mathcal{C}^\infty_c(\mathcal{M})$ and compute
\[
 \begin{array}{ccl}
  \displaystyle \int_{\mathcal{M}}(L_gv)\varphi d\hbox{vol}_g & 
  = & \displaystyle \int_{\mathcal{M}}vL_g\varphi d\hbox{vol}_g
  = \int_{\mathcal{M}}\left(\int_s^\infty (\gamma_{s_1}f)ds_1\right)L_g\varphi d\hbox{vol}_g.
 \end{array}
\]
By Fubini's theorem
\[
\int_{\mathcal{M}}(L_gv)\varphi d\hbox{vol}_g
= \int_s^\infty ds_1\int_{\mathcal{M}}(\gamma_{s_1}f)L_g\varphi d\hbox{vol}_g
= \int_s^\infty ds_1\int_{x\succ \sigma_{s_1}}(\gamma_{s_1}f)L_g\varphi d\hbox{vol}_g.
\]
Using the identity $\psi L_g\varphi - \varphi L_g \psi =
\psi\square_g\varphi - \varphi \square_g \psi =
 \nabla_\mu\left(g^{\mu\nu}(\psi \partial_\nu \varphi - \varphi\partial_\nu \psi)\right)$
 for $\psi = \gamma_{s_1}f$
and Stokes' theorem we find (taking into account the fact that
$\partial \{x\succ \sigma_{s_1}\} = - \sigma_{s_1}$)
\[
 \begin{array}{ccl}
  \displaystyle 
 \int_{\mathcal{M}}(L_gv)\varphi d\hbox{vol}_g & = &\displaystyle 
 \int_s^\infty ds_1\int_{x\succ \sigma_{s_1}}\varphi L_g (\gamma_{s_1}f) d\hbox{vol}_g\\
 & & \displaystyle -
  \int_s^\infty ds_1\int_{\sigma_{s_1}}\langle N,(\gamma_{s_1}f) \nabla \varphi -
  \varphi\nabla (\gamma_{s_1}f)\rangle_g d\mu_g\\
  & = & \displaystyle 0 + \int_s^\infty ds_1\int_{\sigma_{s_1}}\varphi\langle N,
  \nabla (\gamma_{s_1}f)\rangle_g d\mu_g = \int_s^\infty ds_1\int_{\sigma_{s_1}}\varphi\lambda_{s_1}fd\mu_g.
 \end{array}
\]
Hence using the definition of $\lambda_s$ and viewing $d\mu_g$ as a $n$-form, we deduce that
\[
 \begin{array}{ccl}
  \displaystyle 
 \int_{\mathcal{M}}(L_gv)\varphi d\hbox{vol}_g & = &\displaystyle 
 \int_s^\infty ds_1\int_{\underline{\sigma}}\langle N(F_{s_1}),\frac{\partial F}{\partial s}(s_1,\cdot)\rangle_g 
 F_{s_1}^*(\varphi fd\mu_g)\\
 & = &\displaystyle \int_s^\infty\int_{\underline{\sigma}}
 F^*\left(\varphi f \langle N,\cdot \rangle_g \wedge d\mu_g\right)
 \end{array}
\]
But since on $\sigma_{s_1}$, $\langle N,\cdot \rangle_g \wedge d\mu_g = d\hbox{vol}_g$
(again viewing $d\hbox{vol}_g$ as a $(n+1)$-form),
\[
 \int_{\mathcal{M}}(L_gv)\varphi d\hbox{vol}_g
 = \int_s^\infty\int_{\underline{\sigma}} F^*\left(\varphi fd\hbox{vol}_g\right)
 = \int_{x\succ \sigma_s}\varphi fd\hbox{vol}_g.
\]
which proves $L_gv = fY_{\sigma_s}$ in the distribution sense. Since we have obviously
$v = 0$, for $\{x\prec \sigma_s\}$, we deduce $v = \Gamma_sf$ by uniqueness.
Hence (\ref{B=intbeta}) follows.\hfill $\square$

\subsection{Formulation of the dynamics}

We show here a result analogous to Theorem \ref{magic-dynamics} for the nonlinear cubic
Klein--Gordon equation
\begin{equation}\label{u3-kg-equation-on-curved}
 \square_g u + u^3 = 0,
\end{equation}
on a 4-dimensional space-time $\mathcal{M}$ satisfying the hypotheses of Theorem \ref{existence-on-curved},
involving a smooth family of admissible Cauchy hypersurfaces $(\sigma_s)_{s\in \R}$.
We need technical assumptions on $(\sigma_s)_{s\in \R}$, namely:
\begin{equation}\label{hypo1-technique-on-curved}
\exists C_1>0,\quad \forall s\in \R,\forall x\in \sigma_s,\quad 
|\lambda_s(x)| \leq C_1
\end{equation}
and
\begin{equation}\label{hypo2-technique-on-curved}
\exists C_2>0,\quad 
\forall s\in \R,\forall u\in H^1(\sigma_s),\quad \|u\|_{L^6(\sigma_s)}\leq C_2 \|\nabla u\|_{L^2}.
\end{equation}
Note that (\ref{hypo2-technique-on-curved}) is the assumption that the Sobolev embedding $H^1_0(\R^3)\subset L^6(\R^3)$
can be extended on each 3-dimensional manifold $\sigma_s$ uniformly in $s$. This is true if e.g. the Ricci curvature of all 
$\sigma_s$ is uniformly bounded from below and the volumes of all unit balls in $\sigma_s$ are uniformly bounded from below
(see \cite{hebey}).
\begin{theo}\label{theo-magic-dynamics-curved}
Let $(\mathcal{M},g)$ be a 4-dimensional pseudo-Riemannian manifold and $\tau,\rho\in \mathcal{C}^\infty(\mathcal{M})$
satisfying (\ref{hypo-tau-rho-1})--(\ref{further-hypo}).
Let $(\sigma_s)_{s\in \R}$ be a $\underline{\sigma}$-family of admissible Cauchy hypersurfaces which satisfies
(\ref{hypo1-technique-on-curved}) and (\ref{hypo2-technique-on-curved}) and s.t. $\sup_sC(\sigma_s)<+\infty$.
Consider the non autonomous vector field
$V:\R\times \mathcal{E}^1_0(\Sigma_\tau)\longrightarrow  \mathcal{E}^1_0(\Sigma_\tau)$
defined by $V(s,\varphi):= \Phi_{\sigma_s}(0,\lambda_s \varphi^3|_{\sigma_s})$.

Let $I=[t_1,t_2]$ and $J$ be intervals of $\R$ s.t. $\Sigma_{t_1} \prec \sigma_s \prec \Sigma_{t_2}$, $\forall s\in J$
and $u\in \mathcal{F}^1_I(\Sigma_\tau)$.
If $u$ is a solution of (\ref{u3-kg-equation-on-curved}), then $\Theta_{\sigma_s}u$ is a $\mathcal{C}^1$
function of $s\in J$ and satisfies:
\begin{equation}\label{equamagic-dynamics}
 \frac{d(\Theta_{\sigma_s}u)}{ds} + V(s,\Theta_{\sigma_s}u) = 0,\quad \forall s\in J.
\end{equation}
\end{theo}
\emph{Proof} --- First note that $V$ exists and is continuous because of
Corollary \ref{coro-continuity-of-Phi-curved} and of (\ref{hypo1-technique-on-curved}) and
(\ref{hypo2-technique-on-curved}), which imply in particular: $\forall \varphi\in \mathcal{E}^1_{0}(\Sigma_\tau)$, $\forall s\in \R$,
$\lambda_s\varphi^3|_{\sigma_s}\in L^2(\sigma_s)$.
Second let $u\in \mathcal{F}^1_I(\Sigma_\tau)$ and assume that $u$ is a solution of
(\ref{u3-kg-equation-on-curved}).\\
\emph{Step 1} --- We show that $[s\longmapsto \Theta_{\sigma_s}u]$ is continuous, i.e. $\forall s\in J$,
\[
\lim_{s'\rightarrow s}\left(\sup_{t\in I}
E_{(\Theta_{\sigma_{s'}}u)-(\Theta_{\sigma_s}u)}(\Sigma_t)^{1/2}\right) = 0.
\]
Since $(\Theta_{\sigma_{s'}}u)-(\Theta_{\sigma_s}u)\in \mathcal{E}^1_0(\Sigma_t)$, it suffices to prove 
$\lim_{s'\rightarrow s}E_{(\Theta_{\sigma_{s'}}u)-(\Theta_{\sigma_s}u)}(\sigma_{s'}) = 0$
and to apply Corollary \ref{coro-continuity-of-Phi-curved} with $\sigma_{s'}$. But actually
$[\Theta_{\sigma_{s'}}u]_{\sigma_{s'}} = [u]_{\sigma_{s'}}$ so that
$E_{(\Theta_{\sigma_{s'}}u)-(\Theta_{\sigma_s}u)}(\sigma_{s'}) =
E_{u-(\Theta_{\sigma_s}u)}(\sigma_{s'})$. Now observe that $[u-(\Theta_{\sigma_s}u)]_{\sigma_s}=0$
or equivalentely $E_{u-(\Theta_{\sigma_s}u)}(\sigma_s) = 0$.
Thus in particular the result is straightforward in the case where $u$ is smooth with compact spatial support.
The general case follows by proving the existence of a sequence of smooth functions with compact
spatial support which converges to $u$ in the $\mathcal{F}^1_I(\Sigma_\tau)$ topology.
For that purpose first approach $-u^3$ by a sequence of smooth maps with compact spatial
support $(f_\varepsilon)_{\varepsilon>0}$ in $L^1(I,L^2_\ell(\Sigma_\tau))$ and,
for some Cauchy hypersurface $\Sigma$, approach  
$[u]_\Sigma$ by a sequence $(v_\varepsilon,w_\varepsilon)_{\varepsilon>0}$ of smooth maps with compact support in the
$H^1_m(\Sigma)\times L^2(\Sigma)$ topology. For any $\varepsilon>0$ consider the solution $u_\varepsilon$
of $L_gu_\varepsilon = f_\varepsilon$, with the Cauchy data
$[u_\varepsilon]_\Sigma = (v_\varepsilon,w_\varepsilon)$. Then $u_\varepsilon$ is smooth with compact
spatial support and converges to $u$ in $\mathcal{F}^1_I(\Sigma_\tau)$, when $\varepsilon\rightarrow 0$,
because of (\ref{main-estimate-on-curved}).\\
\emph{Step 2} --- We use the generalized Duhamel formula.
First by applying Lemma \ref{lemma-de-Sigma0-vers-tildeSigma} to $u$ and for $\widehat{\Sigma} = \sigma_s$,
we deduce that $s\longmapsto \|u|_{\sigma_s}\|_{H^1_0}$ is bounded.
Hence again because of (\ref{hypo1-technique-on-curved}) and (\ref{hypo2-technique-on-curved}), 
$s\longmapsto \|\lambda_su^3|_{\sigma_s}\|_{L^2}$ is bounded. Thus we can
apply Proposition \ref{duhamel-on-curved}. Then (\ref{duhamel-abstrait}) reads
\[
 u + \Gamma_s(u^3) = \Theta_{\sigma_s}u\quad \hbox{on } \{x\succ \sigma_s\}.
\]
Comparing this identity for two different value $s_1,s_2$ of $s$, we get
\begin{equation}\label{bellesoustraction}
 \Theta_{\sigma_{s_2}}u - \Theta_{\sigma_{s_1}}u = \Gamma_{s_2}(u^3) - \Gamma_{s_1}(u^3)
 \quad \hbox{on } \{x\succ \sigma_{s_1}\}\cap \{x\succ \sigma_{s_2}\}.
\end{equation}
However the r.h.s. of (\ref{bellesoustraction}) can be written by using (\ref{B=intbeta})
\[
 \Gamma_{s_2}(u^3) - \Gamma_{s_1}(u^3) = - \int_{s_1}^{s_2}(\gamma_su^3)ds.
\]
Moreover, since $(\Theta_{\sigma_s}u)|_{\sigma_s} = u|_{\sigma_s}$,
\[
\gamma_su^3 = \Phi_{\sigma_s}(0,\lambda_su^3) = V(s,\Theta_{\sigma_s}u)
\quad \hbox{on } \{x\succ \sigma_s\}
\]
Hence (\ref{bellesoustraction}) implies that the following identity holds on $\{x\succ \sigma_{s_1}\}\cap \{x\succ \sigma_{s_2}\}$:
\begin{equation}\label{belleidentite}
 \Theta_{\sigma_{s_2}}u - \Theta_{\sigma_{s_1}}u 
 + \int_{s_1}^{s_2}V(s,\Theta_{\sigma_s}u)ds =0.
\end{equation}
But since the l.h.s. of (\ref{belleidentite}) is a solution of $L_g\varphi = 0$ on $\mathcal{M}$,
(\ref{belleidentite}) holds actually everywhere on $\mathcal{M}$, by uniqueness. From (\ref{belleidentite}), the
result of the first step and Corollary \ref{coro-continuity-of-Phi-curved}
we then deduce easily (\ref{equamagic-dynamics}).\hfill $\square$

\section{The space of analytic functions over a Banach space}\label{sectionsetting}
\subsection{Analytic functions over a Banach space}\label{sub:functions_spaces} 

Recall that, if $\mathbb{X}$ and $\mathbb{\mathbb{Y}}$ are Banach spaces and
$r\in (0,+\infty)$, $\F_r(\mathbb{X},\mathbb{Y})$ is the space of formal series
$f = \sum_{p=0}^\infty f^{(p)}$ s.t.
$\lc f \rf (r) < +\infty$, where $\lc f \rf (z)$ is given by (\ref{formule-qui-donne-VexX}).
Note that $\left( \F_r(\mathbb{X},\mathbb{Y}),\lc \cdot \rf (r)\right)$ is a Banach space.

Beside the definition of $\lc f\rf$ given by (\ref{formule-qui-donne-VexX}),
we also set, for $k\in \N$,
\[
\lc f\rf ^{(k)}(r):= {d^k\over dz^k}\lc f\rf (z)|_{z=r} \quad
\hbox{and} \quad  \F_r^{(k)}(\mathbb{X},\mathbb{Y}):= \{f\in \F_r(\mathbb{X},\mathbb{Y})|\
\lc f\rf ^{(k)}(r) < +\infty\},
\]
so that for instance $\lc f\rf ^{(1)}(r) = \sum_{p=1}^\infty
p\| f^{(p)}\|_{\otimes} r^{p-1}$. We set
$\F_\infty(\mathbb{X},\mathbb{Y}):= \cap_{r>0}\F_r(\mathbb{X},\mathbb{Y})$ and
$\F_{pol}(\mathbb{X},\mathbb{Y}):= \{f = \sum_{p=0}^Nf^{(p)} |\ N\in \N,
f^{(p)}\in \mathcal{Q}^p\left(\mathbb{X},\mathbb{Y}\right)\}$. Note that we have the dense inclusions
\[
\forall r,R\in (0,\infty), \hbox{ s.t. }r<R, \forall k,\ell\in \N\hbox{ s.t. }k<\ell,\quad
\F_{pol} \subsetneq \F_\infty \subsetneq \F_R \subsetneq \F^{(\ell)}_r\subsetneq \F^{(k)}_r
\subset \F_r.
\]
Indeed if $0<r<R$ and $k\in \mathbb{N}$, we have: $\forall f\in \F_R$,
\begin{equation}\label{comparerN1rNR}
  \lc f\rf ^{(k)}(r) \leq
\Gamma^{(k)}(r,R)\;\lc f\rf (R),\quad\hbox{where }
\Gamma^{(k)}(r,R):= \frac{1}{r^k}\sup_{p\geq k}\frac{p!}{(p-k)!}\left(\frac{r}{R}\right)^p<+\infty.
\end{equation}
In the following we set $\mathbb{Y} =\R$ and:
\begin{defi}
For any $r_0\in (0,\infty]$ and any $k,\ell\in \N$ a \textbf{continuous operator $\T$ from $\F_{(0,r_0)}^{(k)}(\mathbb{X})$ to
$\F_{(0,r_0)}^{(\ell)}(\mathbb{X})$} is a family $\left( \T_r\right)_{0<r<r_0}$, s.t., for any $r\in (0,r_0)$,
$\T_r:\F_r^{(k)}(\mathbb{X})\longrightarrow \F_r^{(\ell)}(\mathbb{X})$
is a continuous linear operator with norm $||\T_r||$ and s.t.,
$\forall r,r'\in (0,r_0)$, if $r<r'$, then the restriction of $\T_r$ to $\F_{r'}^{(k)}(\mathbb{X})$ coincides with $\T_{r'}$.\\
For simplicity we systematically denote each operator $\T_r$ by $\T$ in the following..
\end{defi}

\subsection{Analytic vector fields over $\mathbb{X}$}

\begin{defi}
Elements of $\F_r(\mathbb{X},\mathbb{X})$ are called \textbf{analytic vector fields on $\mathbb{X}$}.
For any $V\in \F_r(\mathbb{X},\mathbb{X})$, we denote by $V\cdot$ the linear
operator acting on $\F_r(\mathbb{X})$ defined by
\[
\forall f\in \F_r(\mathbb{X}), \forall \varphi \in B_\mathbb{X}(r),\quad
\left(V\cdot f\right)(\varphi) = \delta f_\varphi(V(\varphi)),
\]
where, $\forall \varphi \in B_\mathbb{X}(r)$, $\forall \psi\in \mathbb{X}$,
\[
 \delta f_\varphi(\psi):= \lim_{\varepsilon\rightarrow 0}\frac{f(\varphi+\varepsilon\psi)-f(\varphi)}{\varepsilon}
\]
We then set $\lc V\rf \cdot := \lc V\rf (z)\frac{d}{dz}$, a holomorphic vector field on $B_\C(\rho_V)$.
\end{defi}
The previous definition was vague concerning the domain and the target of $V\cdot$.
These points are made more precise by the following result.
\begin{lemm}\label{lemma3.1}
For any $V\in \F_r(\mathbb{X},\mathbb{X})$, the operator $V\cdot$
is continuous from $\F^{(1)}_{(0,r)}(\mathbb{X})$ to $\F_{(0,r)}(\mathbb{X})$ and moreover:
\begin{equation}\label{estop3}
    \forall \rho\in(0,r), \quad \forall f\in\F^{(1)}_\rho(\mathbb{X}), \quad
    \lc V\cdot f\rf (\rho)\leq
    \lc V\rf (\rho) \lc f\rf^{(1)}(\rho)  = (\lc V\rf \cdot \lc f\rf)(\rho).
\end{equation}
\end{lemm}
\emph{Proof} ---
Consider $\rho\in(0,r)$, assume momentaneously that $f \in \F_{pol}(\mathbb{X})$ and
write $f(\varphi) = \sum_{p=0}^Nf^{(p)}(\varphi^{\otimes p})$. Then,
$\forall \varphi\in \mathbb{X}$ such that $||\varphi||_\mathbb{X}\leq r$
we know that $V(\varphi)$ is well defined and, using everywhere the convention $p':= p-1$
and, setting $\varphi_1\otimes \cdots\otimes \varphi_p = \varphi_1\cdots \varphi_p$ for short,
    \[
    \begin{array}{ccl}
   (V\cdot f)(\varphi) & = & \displaystyle
\delta f_\varphi(V(\varphi))
= \sum_{p=1}^N pf^{(p)}(V(\varphi) \underbrace{\varphi\cdots\varphi}_{p'}) \\
     & = & \displaystyle
\sum_{p=1}^N \sum_{q=0}^\infty pf^{(p)} (V^{(q)}(\underbrace{\varphi\cdots\varphi}_q)
\underbrace{\varphi\cdots\varphi}_{p'}) 
= \sum_{m=0}^\infty (V\cdot f)^{(m)}(\underbrace{\varphi\cdots\varphi}_m),
    \end{array}
    \]
where we have set $m=q+p-1 = q+p'$ and, $\forall \varphi_1,\cdots, \varphi_m\in \mathbb{X}$,
    \[
    (V\cdot f)^{(m)}(\varphi_1 \cdots \varphi_m):= \sum_{p=1}^{\sup (N,m+1)}
\frac{1}{m!}\sum_{\sigma\in \mathfrak{S}_m}
pf^{(p)}\left( V^{(m-p')}(\varphi_{\sigma(1)} \cdots \varphi_{\sigma(m-p')})
\varphi_{\sigma(m-p'+1)} \cdots \varphi_{\sigma(m)}\right).
    \]
Hence $|(V\cdot f)^{(m)}(\varphi_1 \cdots \varphi_m)|$ is less than or equal to
(we set $\|\cdot \| = \|\cdot\|_\mathbb{X}$ for shortness):
    \begin{equation}\label{estop4}
\sum_{p=1}^{\sup (N,m+1)} \frac{1}{m!}\sum_{\sigma\in \mathfrak{S}_m}
 p \| f^{(p)}\|_{\otimes} 
\left\|V^{(m-p')}(\varphi_{\sigma(1)}\cdots\varphi_{\sigma(m-p')})
\right\| ||\varphi_{\sigma(m-p'+1)}||\cdots ||\varphi_{\sigma(m)}||
    \end{equation}
and since $\left\|V^{(m-p')}(\varphi_{\sigma(1)}\cdots\varphi_{\sigma(m-p')})\right\|
\leq \|V^{(m-p')}\|_{\otimes} \|\varphi_{\sigma(1)}\|\cdots \|\varphi_{\sigma(m-p')}\|$,
we deduce from the upper bound (\ref{estop4}) that
    \[
    \begin{array}{ccl}
|(V\cdot f)^{(m)}(\varphi_1 \cdots \varphi_m)| & \leq & \displaystyle
\sum_{p=1}^{\sup (N,m+1)} \frac{1}{m!}\sum_{\sigma\in \mathfrak{S}_m}p \|V^{(m-p')}\|_{\otimes} \| f^{(p)}\|_{\otimes}
||\varphi_{\sigma(1)}||\cdots ||\varphi_{\sigma(m)}||\\
 & = & \displaystyle \sum_{p=1}^{\sup (N,m+1)} p \|V^{(m-p')}\|_{\otimes} \| f^{(p)}\|_{\otimes}
||\varphi_1||\cdots ||\varphi_m||.
    \end{array}
    \]
We thus deduce
\begin{equation}\label{newformula-pol}
\|(V\cdot f)^{(m)}\|_{\otimes} \leq \sum_{p=1}^{\sup (N,m+1)}
p \|V^{(m-p')}\|_{\otimes} \| f^{(p)}\|_{\otimes} .
\end{equation}
Hence, by letting $q=m-p'$,
    \[
    \begin{array}{ccl}
    \lc V\cdot f\rf (r) & = & \displaystyle \sum_{m=0}^\infty \|(V\cdot f)^{(m)}\|_{\otimes} r^m \leq
\sum_{m=0}^\infty \sum_{p=1}^{\sup (N,m+1)} p \|V^{(m-p+1)}\|_{\otimes} \| f^{(p)}\|_{\otimes}  r^m\\
     & = & \displaystyle \sum_{q=0}^\infty \sum_{p=1}^N  \|V^{(q)}\|_{\otimes}r^q p\|f^{(p)}\|_{\otimes}
 r^{p-1}
= \displaystyle   \lc V\rf (r) \lc f\rf^{(1)} (r).
    \end{array}
    \]
Thus we obtain (\ref{estop3}) for $f\in \F_{pol}(\mathbb{X})$.
It implies the result by using the density of $\F_{pol}(\mathbb{X})$ in $\F_r^{(1)}(\mathbb{X})$.
\hfill $\square$\\

\noindent
Note that we can extend (\ref{newformula-pol}) \emph{a posteriori} to any
$f\in \F_r^{(1)}(\mathbb{X})$ by density as soon as $\sum_{q=0}^\infty \|V^{(q)}\|r^q<+\infty$,
thanks to (\ref{estop3}). It gives us (still with the convention $p'=p-1$):
\begin{equation}\label{newformula}
\|(V\cdot f)^{(m)}\|_{\otimes} \leq \sum_{p=1}^m  p \|V^{(m-p')}\|_{\otimes} \|f^{(p)}\|_{\otimes} .
\end{equation}
This leads us to the following extension of Lemma \ref{lemma3.1}.
\begin{lemm}\label{reiterer}
Let $k\in \N^*$, $r_0>0$ and $V_1,\cdots ,V_k\in \F_{r_0}(\mathbb{X},\mathbb{X})$.
Then the linear operator $[f\longmapsto V_k\cdot \cdots V_1\cdot \cdot f]$
is continuous from $\F^{(k)}_{(0,r_0)}(\mathbb{X})$ to $\F_{(0,r_0)}(\mathbb{X})$ and $\forall r\in (0,r_0)$,
\begin{equation}\label{NDDF}
\forall f\in \F^{(k)}_r,\quad
\lc V_k\cdot  \cdots V_1\cdot   f\rf (r) \leq
\left(\lc V_k\rf \cdot \cdots \lc V_1\rf \cdot \lc f\rf\right) (r).
\end{equation}
\end{lemm}
\noindent
\emph{Proof} --- For any $a = 1,\cdots ,k$ we write $V_a = \sum_{p=0}^\infty V_a^{(p)}$,
where $\forall p\in \N$, $V_a^{(p)}\in \mathcal{Q}^p(\mathbb{X},\mathbb{X})$. For shortness we set
$X_a^{(p)}:= \| V_a^{(p)}\|_\otimes$, $X_a(z):= \sum_{p\ge 0}X_a^{(p)} z^p$
and $X_a\cdot := X_a(z)\frac{d}{dz}$. We recall that 
$\forall f\in \F, \forall \varphi \in \mathbb{X}$,
$\left(V_a\cdot f\right)(\varphi) = \delta f_\varphi(V_a(\varphi))$.
In the following we assume first that $f\in \F_{pol}$. On the one hand we
observe that, $\forall p\in \N$,
\begin{equation}\label{recursion-to-iterate}
\| (V_k\cdot \cdots V_1\cdot f)^{(p)}\|_{\otimes} \leq
\sum_{p_k'=0}^p\sum_{p_{k-1}'=0}^{p_k}\cdots\sum_{p_1'=0}^{p_2}
p_k\cdots p_1
X_k^{(p-p_k')}X_{k-1}^{(p_k-p_{k-1}')}\cdots X_1^{(p_2-p_1')}\|f^{(p_1)}\|_{\otimes},
\end{equation}
where we systematically denote $p_a':= p_a-1$. This can be proved by recursion on $k$,
by using (\ref{newformula}). On the other hand the coefficients of the decomposition
$(X_k \cdots  X_1\cdot  \lc f\rf )(z) = \sum_{p=0}^\infty
\left(X_k \cdots  X_1\cdot \lc f\rf \right)^{(p)}z^p$
also satisfy similar relations, i.e.
\[
 \left(X_k \cdots  X_1\cdot  \lc f\rf \right)^{(p)}
= \sum_{p_k'=0}^p\sum_{p_{k-1}'=0}^{p_k}\cdots\sum_{p_1'=0}^{p_2}
p_k\cdots p_1
X_k^{(p-p_k')}X_{k-1}^{(p_k-p_{k-1}')}\cdots X_1^{(p_2-p_1')}\| f^{(p_1)}\|_{\otimes},
\]
which can also be proved by a recursion based on the identity
\[
X_a\cdot \left(\sum_{p=0}^\infty A^{(p)}z^p\right) =
\sum_{m=0}^\infty \left(\sum_{p'=0}^mpX_a^{(m-p')}A^{(p)}\right)z^m.
\]
Hence the result follows easily from this identity and (\ref{recursion-to-iterate})
holds for $f\in \F_{pol}(\mathbb{X})$. This can hence be extended to all
$f\in \F_r^{(k)}(\mathbb{X})$ for $r\in (0,r_0)$ by density.
\hfill $\square$.

\section{The time ordered exponential of operators}\label{time_ordered_exponential_of_operators} 
In this section we consider a Lebesgue measurable family $(V_t\cdot)_{t\in I}$ of continuous operators
$V_t\cdot$ from $\F^{(1)}_{(0,r_0)}(\mathbb{X})$ to $\F_{(0,r_0)}(\mathbb{X})$ and we consider the time ordered exponential
\begin{equation}\label{defTexp}
U_{t_1}^{t_2}:= T\hbox{exp}\int_{t_1}^{t_2} d\tau (V_\tau\cdot)
:= \sum_{k=0}^\infty \frac{(V\cdot)^{t_2[k]}_{t_1}}{k!},
\end{equation}
where $(V\cdot)^{t_2[0]}_{t_1}:= 1_{\footnotesize{\hbox{End}(\mathbb{F})}}$ and for $k\geq 1$,
\begin{equation}\label{bdelta(k)t>0}
 (V\cdot)^{t_2[k]}_{t_1}:= k!\int_{{t_1}<\tau_1<\cdots <\tau_k<{t_2}}(V_{\tau_k}\cdots V_{\tau_1}\cdot) d\tau_1\cdots d\tau_k,
\quad \hbox{for }{t_2}>{t_1}
\end{equation}
and
\begin{equation}\label{bdelta(k)t<0}
 (V\cdot)^{t_2[k]}_{t_1}:= (-1)^kk!\int_{{t_2}<\tau_k<\cdots <\tau_1<{t_1}}
(V_{\tau_k}\cdots V_{\tau_1}\cdot) d\tau_1\cdots d\tau_k,
\quad \hbox{for }{t_2}<{t_1}.
\end{equation}
We remark that, for ${t_2}>{t_1}$,
\[
  \frac{(V\cdot)^{t_2[k]}_{t_1}}{k!} = \int_{t_1}^{t_2} d\tau_k\ V_{\tau_k}\cdot
\left(\int_{{t_1}<\tau_1<\cdots <\tau_{k-1}<\tau_k}
(V_{\tau_{k-1}} \cdots V_{\tau_1}\cdot) d\tau_1\cdots d\tau_{k-1}\right)
= \int_{t_1}^{t_2} d\tau\ V_\tau\cdot\frac{(V\cdot)^{\tau[k-1]}_{t_1}}{(k-1)!}.
\]
Hence
\begin{equation}\label{iterative-chrono-int}
 U_{t_1}^{t_2} = 1_{\hbox{\footnotesize{End}}(\mathbb{F})} + \int_{t_1}^{t_2} d\tau\ V_\tau\cdot U_{t_1}^{\tau}.
\end{equation}
A similar reasoning shows that (\ref{iterative-chrono-int}) holds also for ${t_2}<{t_1}$. As a consequence
\begin{equation}\label{fundamental-subtraction}
 U_{t_1}^{t_2+h} - U_{t_1}^{t_2} = \int_{t_1}^{{t_2}+h} d\tau\ V_\tau\cdot U_{t_1}^{\tau} - \int_{t_1}^{t_2} d\tau\ V_\tau\cdot U_{t_1}^{\tau}
= \int_{t_2}^{t_2+h} d\tau\ V_\tau\cdot U_{t_1}^{\tau}.
\end{equation}

\subsection{Existence of $U_{t_1}^{t_2}$}\label{sub-main_result} 
In the following, for any vector field $X$ on $B_\C(r)$, we denote by
$(t,z)\longmapsto e^{-tX}(z)$ the map which is equal to the solution $\gamma$ of
\[
\left\{ \begin{array}{ccl}
\displaystyle{\partial \gamma\over \partial t}(t,z) & = & \displaystyle{-X(\gamma(t,z))}\\
\gamma(0,z) & = & z.
\end{array}\right.
\]
\begin{theo}\label{theo.produit.ordonne.general}
Let $r_0\in (0,+\infty]$ and $I\subset \R$ be an interval.
Let $(V_t)_{t\in I}$ be a normal family of analytic vector fields in $\F_{r_0}(\mathbb{X},\mathbb{X})$
and let $X = \sum_{k=0}^\infty X_kz^k\in \F_{r_0}(\R)$ 
s.t. $\forall t\in I$, $\forall p\in \N$, $0\leq \lc V_t^{(p)}\rf \leq X_p$.
Assume that:
\begin{equation}\label{hypo-mesurable}
\forall r\in (0,r_0), \forall f\in \F_r^{(1)}(\mathbb{X}), \quad [I\ni t\longmapsto V_t\cdot f\in \F_r(\mathbb{X})]
\hbox{ is measurable}.
\end{equation}
Let $R\in (0,r_0)$. Then $\forall t_1,t_2\in I$ s.t. $e^{-|t_2-t_1|X}(R)$ exists and is positive,
the operator $U_{t_1}^{t_2}:= T\exp\left( \int_{t_1}^{t_2} d\tau V_\tau\cdot \right)$
defined by (\ref{defTexp}) is a bounded operator from $\F_R(\mathbb{X})$ to $\F_{e^{-|t_2-t_1|X}(R)}(\mathbb{X})$ with 
a norm less than $1$ i.e.
\begin{equation}\label{maininequa}
\forall t\in [t_1,t_2],
\forall f\in  \F_R(\mathbb{X}), \quad \lc U_{t_1}^{t_2}f \rf (e^{-|t_2-t_1|X}(R))
\leq \lc f\rf (R).
\end{equation}
Moreover for any $\overline{R}$ s.t. $R<\overline{R}<r_0$ and $f\in \F_{\overline{R}}(\mathbb{X})$,
the map $t\longmapsto U_{t_1}^{t}f$ is locally Lipschitz continuous from $[t_1,t_2]$ to $\F_{e^{-|t_2-t_1|X}(R)}(\mathbb{X})$. 
\end{theo}
\emph{Proof of theorem \ref{theo.produit.ordonne.general}} ---
W.l.g. we assume throughout the proof that $t_1=0<T=t_2$ and study $U_0^t$ for $0<t\leq T$.
The proof is divided in several steps which follow.\\
\noindent
\emph{Step 1} --- For $r\in (0,r_0)$, $k\in \N$, $f\in \F(\mathbb{X})$ and $t\in I$ we estimate the norm
in $\F_r(\mathbb{X})$ of $(V\cdot)^{t[k]}_0 f$. For $t\in [0,T]$,
we start from Expression (\ref{bdelta(k)t>0}) for $(V\cdot)^{t[k]}_0 f$
and we use Lemma \ref{reiterer} with $V_a\cdot  = V_{\tau_a}\cdot$.
This gives us
\[
 \begin{array}{ccl}
\displaystyle \lc (V\cdot)^{t[k]}_0 f\rf(r) & \leq &
\displaystyle k!\int_{0<\tau_1<\cdots <\tau_k<t}
\lc V_{\tau_k}\cdots V_{\tau_1}\cdot f\rf(r) d\tau_1\cdots d\tau_k\\
& \leq & \displaystyle k!\int_{0<\tau_1<\cdots <\tau_k<t}
(\lc V_{\tau_k}\rf\cdots \lc V_{\tau_1}\rf\cdot \lc f\rf)(r) d\tau_1\cdots d\tau_k\\
& \leq & \displaystyle  k!\int_{0<\tau_1<\cdots <\tau_k<t} (X\cdot)^k \lc f\rf (r)
d\tau_1\cdots d\tau_k,
 \end{array}
\]
which implies, by using $k!\int_{0<\tau_1<\cdots<\tau_k<t}d\tau_1\cdots d\tau_k = t^k$:
\begin{equation}\label{Nek}
\lc (V\cdot)^{t[k]}_0 f\rf(r) \leq t^k (X\cdot)^k \lc f\rf (r).
\end{equation}
Hence we see how to derive a sufficient condition for the series
$U_0^tf = \sum_{k=0}^\infty \frac{1}{k!}(V\cdot)^{t[k]}_0 f$ to be convergent
in some space $\F_r$: it suffices to find some $r$ which satisfies
\begin{equation}\label{sufficient-cond-conv}
\sum_{k=0}^\infty {t^k\over k!} (X\cdot)^k \lc f\rf(r) < +\infty.
\end{equation}
Then this implies by (\ref{Nek}) that $\sum_{k=0}^\infty \frac{1}{k!}
\lc(V\cdot)^{t[k]}_0 f \rf(r) < +\infty$ and hence the existence of
$U_0^tf$.\\

\noindent
\emph{Step 2} --- We show that, if $R\in (0,r_0)$, $T>0$ and $e^{-TX}(R)>0$, condition (\ref{sufficient-cond-conv})
is satisfied with $t=T$ and $r = e^{-TX}(R)$. Actually we will show that 
\begin{equation}\label{sufficient-cond-holds}
\forall f\in \F(\mathbb{X}), \quad
\sum_{k=0}^\infty {T^k\over k!} ((X\cdot)^k \lc f\rf)\left(e^{-TX}(R)\right) = \lc f\rf(R).
\end{equation}
For that purpose we use the following lemma, the proof of which is given below.
In the following, for $t,r>0$, we set $\overline{B}_\C(r):= \{z\in \C|\ |z|\leq r\}$.

\begin{lemm}\label{champdevecteur}\sl{
Let $X:B_\C(r_0)\longmapsto \C$ be an holomorphic vector field different from 0. Assume that
\[
X(z) = \sum_{k=0}^\infty X_kz^k,\quad \hbox{where }X_k\geq 0, \forall k\in \N.
\]
Let $\rho\in(0,r_0)$ and $T>0$ such that $e^{TX}(\rho)$ exists. Then the flow map
\[
\begin{array}{ccc}
\overline{B}_\C(T)\times \overline{B}_\C(\rho) & \longrightarrow & \C\\
(\tau,z) & \longmapsto & e^{\tau X}(z)
\end{array}
\]
is well defined  and holomorphic and in particular
\begin{equation}\label{controleflot}
\forall (\tau,z)\in \overline{B}_\C(T)\times \overline{B}_\C(\rho),\quad
|e^{\tau X}(z)| \leq e^{|\tau|X}(|z|)\leq e^{TX}(\rho).
\end{equation}
}
\end{lemm}
\noindent Consider any $R\in (0,r_0)$, $0<t\leq T$ s.t. $e^{-tX}(R)>0$: then 
$e^{tX}\left(e^{-tX}(R)\right)$ exists since it is nothing but $R$.
Hence we can apply Lemma \ref{champdevecteur} with $\rho = e^{-TX}(R)$. It implies in particular that, for any
holomorphic function $H$ on $\overline{B}_\C(R) = \overline{B}_\C(e^{TX}(\rho))$, the map
\[
\begin{array}{ccc}
\overline{B}_\C(T)\times \overline{B}_\C(e^{-TX}(R)) & \longrightarrow & \C\\
(\tau,z) & \longmapsto & H\left( e^{\tau X}(z) \right)
\end{array}
\]
is well defined and is analytic. Hence the following expansion holds:
\begin{equation}\label{taylor}
\forall (\tau,z) \in \overline{B}_\C(T)\times \overline{B}_\C(e^{-TX}(R)),\quad
H\left( e^{\tau X}(z) \right) =
\sum_{k=0}^\infty \left.\frac{d^kH\left( e^{s X}(z) \right)}{(ds)^k}\right|_{s=0}
 {\tau^k\over k!},
\end{equation}
the series on the r.h.s. being absolutely convergent for any $\tau\in \overline{B}_\C(T)$.
However because of the identity
$\left(\frac{d}{ds}\right)^k\left[H\left( e^{sX}(z) \right)\right] = ((X\cdot)^k H)(e^{sX}(z))$,
which can be proved by recursion over $k$, we deduce from (\ref{taylor}) that
\[
\forall (\tau,z) \in \overline{B}_\C(T)\times \overline{B}_\C(e^{-TX}(R)),\quad
H\left( e^{\tau X}(z) \right) = \sum_{k=0}^\infty ((X\cdot)^k H)(z) {\tau^k\over k!}.
\]
By specializing this relation to $(\tau,z) = (T,e^{-TX}(R))$ we deduce that the power series
$\sum_{k=0}^\infty {T^k\over k!}((X\cdot)^k h)(e^{-TX}(R))$ is absolutely convergent and satisfies the identity
\begin{equation}\label{identity=}
H(R) = \sum_{k=0}^\infty {T^k\over k!}\left((X\cdot)^k H\right)(e^{-TX}(R)).
\end{equation}
Hence by using (\ref{identity=}) with $H(z) = \lc f\rf(z)$ we obtain (\ref{sufficient-cond-holds}).
This shows that the series $T \hbox{exp}\left( \int_0^T d\tau V_\tau\cdot\right) f$
converges in $\F_{e^{-TX}(R)}(\mathbb{X})$. Moreover we deduce using (\ref{Nek}) and (\ref{sufficient-cond-holds}) the
following estimate:
\begin{equation}\label{provisory-conclusion}
\left[ T \hbox{exp}\left( \int_0^T d\tau V_\tau\cdot\right) f \right](e^{-TX}(R))
\leq \sum_{k=0}^\infty  \frac{1}{k!}\left[(V\cdot)^{T[k]}_0 f \right]\left(e^{-TX}(R)\right)
\leq  \lc f\rf(R).
\end{equation}
Lastly we remark that hypothesis $e^{-TX}(R)>0$ obviously implies $e^{-tX}(R)>0$, $\forall t\in [0,T]$
so that Conclusion (\ref{provisory-conclusion}) holds also if we replace $T$ by $t\in [0,T]$.
This implies (\ref{maininequa}).\\

\noindent
\emph{Step 3 } --- Let us prove the local Lipschitz continuity of $t\longmapsto U_0^tf$,
for $f\in \F_{\overline{R}}(\mathbb{X})$, where $R<\overline{R}<r_0$.
Let $t\in[0,T]$ and $h\in \mathbb{R}$ s.t. $t+h\in [0,T]$. Then it follows from
(\ref{fundamental-subtraction}) and (\ref{estop3}) that
\[
\begin{array}{ccl}
\displaystyle \lc (U_0^{t+h}-U_0^t)f\rf\left(e^{-TX}(R)\right) & \leq & \displaystyle 
\int_t^{t+h}d\tau \lc V_\tau\cdot U_0^{\tau}f\rf\left(e^{-TX}(R)\right) \\
& \leq & \displaystyle |h|X\left(e^{-TX}(R)\right)\sup_{t<\tau<t+h}
\lc U_0^{\tau}f\rf^{(1)} \left(e^{-TX}(R)\right).
\end{array}
\]
However by observing that $e^{-tX}(R)<e^{-tX}(\overline{R})$ because $R<\overline{R}$ we
deduce from (\ref{comparerN1rNR}) that 
$\lc g\rf^{(1)}\left(e^{-TX}(R)\right) \leq \Gamma^{(1)}(R,\overline{R})\,
\lc g\rf\left(e^{-TX}(\overline{R})\right)$,
$\forall g\in \F_{e^{-TX}(\overline{R})}(\mathbb{X})$. Applying this for
$g = U^\tau_0f$,
\[
\lc (U_0^{t+h}-U_0^t)f\rf\left(e^{-TX}(R)\right) \leq 
|h|X\left(e^{-TX}(R)\right)\Gamma^{(1)}(R,\overline{R})
\sup_{t<\tau<t+h} \lc U_0^{\tau}f\rf\left(e^{-TX}(\overline{R})\right)
\]
and by using (\ref{maininequa}) with $\overline{R}$ instead of $R$:
\[
\lc (U_0^{t+h}-U_0^t)f\rf\left(e^{-TX}(R)\right) \leq 
|h|X\left(e^{-TX}(R)\right)\Gamma^{(1)}(R,\overline{R})\; \lc f\rf\left(\overline{R}\right).
\]
\hfill $\square$\\

\noindent
\emph{Proof of lemma \ref{champdevecteur}} ---
We first show that $(\tau,z) \longmapsto e^{\tau X}(z)$ is defined and satisfies (\ref{controleflot})
over $B_\C(T)\times B_\C(\rho)$. Fix some $z\in B_\C(\rho)$ and $\tau\in B_\C(T)$.
Then $\exists \varepsilon_0>0$ s.t. $\forall \varepsilon\in (0,\varepsilon_0]$,
\[
|z| \leq \rho - \varepsilon \quad \hbox{and}\quad
|\tau|\leq T_\varepsilon:= {T\over 1+\varepsilon}.
\]
We also let $\lambda\in S^1\subset \C$ such that $\tau = |\tau|\lambda$, where $0< |\tau| \leq T_\varepsilon$.
We introduce the notations:
\[
\left\{ \begin{array}{ccll}
f_\varepsilon(t) &:= & e^{t(1+\varepsilon)X}(|z|+\varepsilon) & \forall t\in [0,T_\varepsilon]\\
\gamma(t) &:= & e^{t\lambda X}(z) & \forall t\in [0, \overline{t}) \\
g(t) & := & |\gamma(t)| & \forall t\in [0, \overline{t})
\end{array}\right.
\]
where $\overline{t}$ is the positive maximal existence time for $\gamma$. Note that $f_\varepsilon$
is defined on $[0,T_\varepsilon]$ because of the assumption that $e^{TX}(R)$ exists. Our first task is to show that the set:
\[
A_\varepsilon:= \{t\in [0,T_\varepsilon]\cap [0,\overline{t}) |\ g(t) - f_\varepsilon(t) \geq 0\}
\]
is actually empty. Let us prove it by contradiction and assume that
$A_\varepsilon\neq \emptyset$. Then there exists $t_0:= \inf
A_\varepsilon$. Note that $g(0) - f_\varepsilon(0) = - \varepsilon
<0$, hence we deduce from the continuity of $g- f_\varepsilon$ that
$t_0\neq 0$ and $g(t_0) = f_\varepsilon(t_0)$. Moreover since
$f_\varepsilon(0) = \varepsilon$ and $f_\varepsilon$ is increasing
because $X(r)>0$ for $r>0$ we certainly have $g(t_0) =
f_\varepsilon(t_0) >0$. We now observe that
\[
\forall z\in \C^*,\quad {\langle \lambda X(z),z\rangle \over |z|}
= \left\langle \lambda \sum_{k=0}^\infty X_kz^k,{z\over |z|}\right\rangle \leq \sum_{k=0}^\infty X_k |z|^k = X(|z|).
\]
Hence for all $t\geq 0$ s.t. $g(t)\neq 0$,
\[
g'(t) = {\langle \lambda X(\gamma(t)),\gamma(t)\rangle \over |\gamma(t)|} \leq X(|\gamma(t)|) = X(g(t))
\]
and hence in particular, since $g(t_0) \neq 0$,
\[
g'(t_0) \leq X(g(t_0)) = X(f_\varepsilon(t_0)) = {f_\varepsilon'(t_0) \over 1+\varepsilon} < f_\varepsilon'(t_0).
\]
Thus since $f'_\varepsilon - g' $ is continuous $\exists t_1\in (0,t_0)$
s.t. $\forall t\in [t_1,t_0]$, $f'_\varepsilon(t) - g'(t)\geq 0$. Integrating this inequality over $[t_1,t_0]$ we obtain
\[
g(t_1) - f_\varepsilon(t_1) = \left( f_\varepsilon(t_0) - g(t_0)\right) - \left( f_\varepsilon(t_1) - g(t_1)\right)
= \int_{t_1}^{t_0} \left(f'_\varepsilon(t) - g'(t)\right) dt \geq 0,
\]
i.e. $t_1\in A_\varepsilon$, a contradiction. \\

\noindent
Hence $A_\varepsilon = \emptyset$. Note that this implies automatically that $\overline{t}> T_\varepsilon$.
Indeed if we had $\overline{t}\leq T_\varepsilon$ this would imply that $g$ is not bounded in
$[0,\overline{t}) \subset [0,T_\varepsilon]$, but since $f_\varepsilon$ is bounded on $[0,T_\varepsilon]$
we could then find some time $t\in [0,\overline{t})$ s.t. $g(t)\geq f_\varepsilon(t)$,
which would contradict the fact that $A_\varepsilon = \emptyset$.
Thus we deduce that $\forall t\in [0,T_\varepsilon]$, $g(t) < f_\varepsilon(t)$, i.e.
\[
\forall t\in [0,T_\varepsilon],\quad |e^{\lambda tX}(z)| < e^{(1+\varepsilon)tX}(|z|+\varepsilon).
\]
In other words for all $\tau = \lambda t\in B_\C(T)$ and all $z\in B_\C(\rho)$ we found that
$\forall \varepsilon\in (0,\varepsilon_0]$, $|e^{\tau X}(z)|\leq e^{(1+\varepsilon)|\tau|X}(|z|+\varepsilon)$.
Letting $\varepsilon$ goes to 0, we deduce the estimate (\ref{controleflot}) for
$(\tau,z)\in B_\C(T)\times B_\C(\rho)$.
Lastly this estimate forbids the flow to blow up on $\overline{B}_\C(T)\times \overline{B}_\C(\rho)$.
Hence the result and (\ref{controleflot}) can be extended to this domain by continuity.
\hfill $\square$

\section{Proof of the Main Theorem}

We first prove the following strengthening of Theorem \ref{theo.produit.ordonne.general}
(with stronger hypotheses).
\begin{theo}\label{theo.produit.ordonne.general.C1}
Let $r_0\in (0,+\infty]$ and $I\subset \R$ be an interval.
Let $(V_t)_{t\in I}$ be a normal family of analytic vector fields in $\F_{r_0}(\mathbb{X},\mathbb{X})$
and let $X = \sum_{k=0}^\infty X_kz^k\in \F_{r_0}(\R)$ 
s.t. $\forall t\in I$, $\forall p\in \N$, $0\leq \lc V_t^{(p)}\rf \leq X_p$.
Assume that:
\begin{equation}\label{hypo-plusfort-continu}
I\times B_\mathbb{X}(r_0)\ni (t,\varphi)\longmapsto V_t(\varphi)\in \mathbb{X}
\hbox{ is continuous}.
\end{equation}
Let $R,\overline{R}\in \mathbb{R}$ s.t. $0<R<\overline{R} <r_0$.
Let $f\in \mathbb{F}_{\overline{R}}(\mathbb{X})$.
Let $t_1,t_2\in I$ s.t. $e^{-|t_2-t_1|X}(R)>0$ and
let $\varphi\in \mathcal{C}^1([t_1,t_2],\mathbb{X})$ s.t.
$\|\varphi(t)\|_\mathbb{X}\leq e^{-|t-t_1|X}(R)$, $\forall t\in [t_1,t_2]$.
Then the map
\[
 \begin{array}{ccl}
  [t_1,t_2] & \longrightarrow & \mathbb{X}\\
t & \longmapsto & \left(U_{t_1}^tf\right)(\varphi(t))
 \end{array}
\]
is $\mathcal{C}^1$ and satisfies
\begin{equation}\label{StF-est-C1}
 \frac{d}{dt}\left((U_{t_1}^tf)(\varphi(t))\right) = (V_t\cdot U_{t_1}^tf)(\varphi(t))
+ \delta(U_{t_1}^tf)_{\varphi(t)}\left(\frac{d\varphi(t)}{dt}\right).
\end{equation}
\end{theo}
\emph{Proof} --- W.l.g. we assume $t_1=0<T=t_2$. 
Let $f\in \F_{\overline{R}}(\mathbb{X})$ and, for $t\in [0,T]$, set $f_t:=U_0^tf$.
By Theorem \ref{theo.produit.ordonne.general}
we know that $f_t\in \F_{e^{-tX}(\overline{R})}(\mathbb{X})$, $\forall t\in [0,T]$.
We first show that, $\forall t\in [0,T]$, $(\tau,\varphi) \longmapsto (V_\tau\cdot f_\tau) (\varphi)$
is continuous on $[0,t]\times B_{\mathbb{X}}(0,e^{-tX}(R))$. For that purpose, for $\tau, \tau + \sigma
\in [0,t]$ and $\varphi,\psi\in B_{\mathbb{X}}(0,e^{-tX}(R))$ we evaluate the difference
$(V_{\tau+\sigma}\cdot f_{\tau+\sigma})(\psi) - (V_\tau\cdot f_\tau)(\varphi)$.
We split this quantity as the sum of three terms:
\[
 ( V_{\tau+\sigma}\cdot f_{\tau+\sigma})(\psi) - (V_\tau\cdot f_\tau)(\varphi)
= \delta (f_{\tau+\sigma})_\psi(V_{\tau+\sigma}(\psi))
- \delta (f_{\tau})_\varphi(V_{\tau}(\varphi))
= \Delta_1 + \Delta_2 + \Delta_3,
\]
where
\[
 \begin{array}{ccl}
  \Delta_1 & := & \delta (f_{\tau+\sigma})_\psi(V_{\tau+\sigma}(\psi) - V_{\tau}(\varphi))\\
  \Delta_2 & := & \left( \delta (f_{\tau+\sigma})_\psi
- \delta (f_{\tau+\sigma})_\varphi\right)(V_{\tau}(\varphi))\\
  \Delta_3 & := & \delta(f_{\tau+\sigma} - f_{\tau})_\varphi(V_{\tau}(\varphi)).
 \end{array}
\]
To evaluate $\Delta_1$ and $\Delta_3$ we will use the following inequality (for all $r>0$):
\begin{equation}\label{simple-consequence-estop3}
\forall g\in \mathbb{F}_r^{(1)}(\mathbb{X}),
\forall \varphi\in B_{\mathbb{X}}(r), \forall Z\in \mathbb{X},\quad |\delta g_\varphi(Z)|
\leq \lc Z\cdot g\rf\left(\|\varphi\|_\mathbb{X}\right)
= \|Z\|_\mathbb{X} \lc g\rf^{(1)}\left(\|\varphi\|_\mathbb{X}\right),
\end{equation}
which follows by applying Lemma \ref{lemma3.1}, (\ref{estop3}) with
$V$ being the constant vector field $[\varphi\longmapsto Z]$.

We note that (\ref{simple-consequence-estop3}), $\|\psi\|_{\mathbb{X}} < e^{-tX}(R)$ and
Inequality (\ref{comparerN1rNR}) imply
\[
\begin{array}{ccl}
 |\Delta_1| & \leq & \|V_{\tau+\sigma}(\psi) - V_{\tau}(\varphi)\|_\mathbb{X}
 \lc f_{\tau+\sigma}\rf^{(1)}\left(\|\psi\|_\mathbb{X}\right)\\
& \leq & \|V_{\tau+\sigma}(\psi) - V_{\tau}(\varphi)\|_\mathbb{X}
 \Gamma^{(1)}(e^{-tX}(R),e^{-tX}(\overline{R}))\lc f_{\tau+\sigma}\rf(e^{-tX}(\overline{R}))
 \end{array}
\]
and hence $\Delta_1$ converges to 0 as $\sigma\rightarrow 0$ and $\|\psi-\varphi\|_\mathbb{X}\rightarrow 0$
because of (\ref{hypo-plusfort-continu}). We decompose and split $\Delta_2$:
\[
\begin{array}{ccl}
 \Delta_2 & := &  \displaystyle \sum_{p=0}^\infty pf^{(p)}_{\tau+\sigma}(V_\tau(\varphi)\underbrace{\psi\cdots \psi}_{p-1}) -
pf^{(p)}_{\tau+\sigma}(V_\tau(\varphi)\underbrace{\varphi\cdots \varphi}_{p-1})\\
& = & \displaystyle \sum_{p=0}^\infty p\sum_{j=1}^{p-1} f^{(p)}_{\tau+\sigma}
(V_\tau(\varphi)\underbrace{\psi\cdots \psi}_{j-1}(\psi-\varphi)\underbrace{\varphi\cdots \varphi}_{p-1-j})
\end{array}
\]
We deduce that, by setting $M:= \sup(\|\varphi\|_\mathbb{X},\|\psi\|_\mathbb{X})$,
\[
 |\Delta_2| \leq \|V_\tau(\varphi)\|_\mathbb{X} \|\psi-\varphi\|_\mathbb{X} \sum_{p=0}^\infty p(p-1)
\| f^{(p)}_{\tau+\sigma}\|_{\otimes} M^{p-2}
= \|V_\tau(\varphi)\|_\mathbb{X} \|\psi-\varphi\|_\mathbb{X} \lc f_{\tau+\sigma}\rf^{(2)}\left(M\right).
\]
Hence by using $M \leq e^{-tX}(R) < e^{-tX}((R+\overline{R})/2)$ and
Inequality (\ref{comparerN1rNR}), we deduce that $\Delta_2$ tends to 0 when $\|\psi-\varphi\|
\rightarrow 0$. Lastly using again (\ref{simple-consequence-estop3}) we have
\[
 |\Delta_3| \leq \|V_{\tau}(\varphi)\|_\mathbb{X} \lc f_{\tau+\sigma} - f_{\tau}\rf^{(1)}\left(\|\varphi\|_\mathbb{X}\right),
\]
which implies also that $\Delta_3$ tends to 0 when $\sigma \rightarrow 0$ by applying
Theorem \ref{theo.produit.ordonne.general} with $(R+\overline{R})/2$ in place of $R$
(since $(R+\overline{R})/2<\overline{R}$ and $f\in \mathbb{F}_{\overline{R}}$ the
map $\tau\longmapsto f_\tau$ is continuous from $[0,t]$ to
$\mathbb{F}_{e^{-tX}(\frac{R+\overline{R}}{2})}(\mathbb{X})$ and hence to
$\mathbb{F}^{(1)}_{e^{-tX}(R)}(\mathbb{X})$ by Inequality (\ref{comparerN1rNR})). 

Hence we conclude that
$(V_{\tau+\sigma}\cdot f_{\tau+\sigma})(\psi) - (V_\tau\cdot f_\tau)(\varphi)$
converges to 0 when $\sigma \rightarrow 0$ and $\|\psi-\varphi\| \rightarrow 0$, which
proves the continuity of  $(\tau,\varphi) \longmapsto (V_\tau\cdot f_\tau) (\varphi)$.

An easy consequence is that the r.h.s. of (\ref{StF-est-C1}) is continuous. Thus it suffices to prove
(\ref{StF-est-C1}) in order to conclude.
Let $h\neq 0$, then using (\ref{fundamental-subtraction}):
\begin{equation}\label{decomposition-Ftphit}
 \begin{array}{ccl}
  \frac{1}{h}\left[f_{t+h}(\varphi(t+h)) - f_t(\varphi(t))\right] & = &
\displaystyle \frac{1}{h}\left[f_{t+h}(\varphi(t+h)) - f_t(\varphi(t+h))\right]
+ \frac{1}{h}\left[f_t(\varphi(t+h)) - f_t(\varphi(t))\right]\\
& = & \displaystyle \frac{1}{h}\int_t^{t+h}d\tau(V_\tau\cdot f_\tau)(\varphi(t+h))
+ \frac{1}{h}\left(f_t(\varphi(t+h)) - f_t(\varphi(t))\right)
 \end{array}
\end{equation}
When $h\rightarrow 0$ the first term in the r.h.s. of (\ref{decomposition-Ftphit})
converges to $(V_t\cdot f_t)(\varphi(t))$
because of the continuity of $(\tau,\varphi) \longmapsto (V_\tau\cdot f_\tau) (\varphi)$.
The second term in the r.h.s. of (\ref{decomposition-Ftphit}) converges to
$\delta(f_t)_{\varphi(t)}\left(\frac{d\varphi(t)}{dt}\right)$ because of
(\ref{seriesC2}). Hence
the r.h.s. of (\ref{decomposition-Ftphit}) converges to
$(V_t\cdot f_t)(\varphi(t)) + \delta(f_t)_{\varphi(t)}\left(\frac{d\varphi(t)}{dt}\right)$
when $h\rightarrow 0$, which proves (\ref{StF-est-C1}). \hfill $\square$\\

\noindent
\emph{Proof of the Theorem \ref{theo-intro}} --- On a flat space-time with a general real analytic
nonlinearity we first use
Proposition \ref{proposition-fund-Section2} which provides us with a normal family
of analytic vector fields $(V_t)_{t\in I}$ satisfying (\ref{hypo-plusfort-continu})
and using Theorem \ref{theo-magic-dynamics} we obtain a $\mathcal{C}^1$ map $\varphi(t) = \Theta_tu$
which satisfies (\ref{magic-dynamics}). We can thus
apply Theorem \ref{theo.produit.ordonne.general.C1} to these data and deduce:
\[
 \frac{d}{dt}\left((U_{t_1}^tf)(\Theta_tu)\right) = (V_t\cdot U_{t_1}^tf)(\Theta_tu)
+ \delta(U_{t_1}^tf)_{\varphi(t)}\left(-V_t(\Theta_tu)\right) =0.
\]
Hence the results follows.

A similar result holds for the Klein--Gordon $\square_gu+u^3=0$
on a 4-dimensional hyperbolic pseudo-Riemannian manifold, by using
Theorem \ref{theo-magic-dynamics-curved} and Theorem \ref{theo.produit.ordonne.general.C1}.
\hfill $\square$

\section{Comparison with quantum field theory}\label{sec-comparison}

The space $\mathbb{F}$ shares some analogies with the Fock spaces used by physicists
in the quantum field theory. In the following we set $N:= \hbox{dim}E$, we let $(e_1,\cdots, e_N)$
be a basis of $E$ and we use the affine coordinates
$E\ni w \longmapsto w^i\in \R$, for $i=1,\cdots, N$, in this basis.
First assume that $s>n/2$, so that $\mathcal{E}_0^s$ embedds continuously in continuous
functions. Then for all $x\in \mathcal{M}$ and $i = 1,\cdots ,N$ we define the continuous linear map
$\bphi^i(x):\mathcal{E}_0^s\longrightarrow \R$
(equivalentely $\bphi^i(x)\in (\mathcal{E}_0^s)^*\subset \mathbb{F}$) by
\[
 \begin{array}{cccc}
  \bphi^i(x): & \mathcal{E}_0^s & \longrightarrow & \R\\
& \varphi & \longmapsto & \varphi^i(x).
 \end{array}
\]
If $s$ is arbitrary we define $\bphi^i$ as a distribution on $\mathcal{M}$,
with values in $(\mathcal{E}_0^s)^*\subset \mathbb{F}$ by
\[
 \begin{array}{cccc}
  \bphi^i: & \mathcal{C}^\infty_c(\mathcal{M}) & \longrightarrow & (\mathcal{E}_0^s)^*\\
& f & \longmapsto & \left[\int_\mathcal{M}f(x)\bphi^i(x)dx : \varphi\longmapsto \int_\mathcal{M}f(x)\varphi^i(x)dx\right]
 \end{array}
\]
Similarly we define $\frac{\partial \bphi^i}{\partial x^\mu}$
as a distribution with values in $(\mathcal{E}_0^s)^*$.
More generally, assuming that $s$ is s.t. we can make sense of $N(\varphi,\partial \varphi)$,
we define the $\F$-valued distribution $N(\bphi,\partial \bphi)$. Note that the constant functional $\textbf{1}$
equal to 1 on $\mathcal{E}_0^s$ plays a role analogous to the vacuum.

As an algebra of functions (on $\mathcal{E}^s_0$) $\F$ acts linearly on itself by multiplication:
to each $g\in \mathbb{F}$ we associate the multiplication linear operator
$[f \longmapsto g f]\in \hbox{End}(\F)$. This defines a natural
embedding $\F\hookrightarrow \hbox{End}(\F)$ and all previous $\F$-valued distributions
$\bphi^i$, $\frac{\partial \bphi^i}{\partial x^\mu}$, $N^i(\bphi,\partial \bphi)$ can
also be viewed as $\hbox{End}(\F)$-valued distributions.

Another important type of $\hbox{End}(\F)$-valued distribution is:
\[
 \begin{array}{cccc}
 \bphi^+_i: & \mathcal{C}^\infty_c(\mathcal{M}) & \longrightarrow & \hbox{End}(\mathbb{F})\\
& f & \longmapsto & \left[\int_\mathcal{M}f(y)\bphi^+_i(y)dy : f\longmapsto
\delta f_{\int_\mathcal{M}f(y)G_ye_idy} \right].
 \end{array}
\]
Here, in the case where $\mathcal{M}=M$ is a flat space-time, $G_y$ is defined by:
$\forall x,y\in M$, $G_y(x):= G(x-y)$, where $G$ is the distribution
defined in (\ref{Gdecomposition}).
In the case where $\mathcal{M}$ is a curved globally hyperbolic space-time and if $L = \square_g$,
$G_y$ is defined in Section \ref{subsec-gen-Duhamel}, i.e.
is the solution of $\square_gG_y + m^2G_y = 0$
with the Cauchy conditions $G_y|_{\sigma} = 0$ and $\langle N,\nabla G_y\rangle_g|_{\sigma} = \delta_y$,
for any Cauchy hypersurface $\sigma$ which contains $y$.
In both case it may be useful to set $G(x,y):= G_x(y)$.

Hence for any $f\in \mathcal{C}^\infty_c(M)$,
$\int_Mf(y)\bphi^+_i(y)dy$ is the analytic first order operator associated with the constant
vector field equal to $\int_Mf(y)G_ye_idy\in \mathcal{E}^s_0$ everywhere.
Intuitively one may think that the notation $\bphi^+_i(y)$
would represent the first order operator $f\longmapsto \delta f_{G_ye_i}$
associated with the constant vector field
$G_ye_i$, if $G_ye_i$ would be in $\mathcal{E}^s_0$ (but it does not here if $s>n/2$).

This language allows us to express the operator $V_t\cdot$ of our Main Theorem as:
\[
V_t\cdot:= \int_{\mathbb{R}^n}d\vec{y}\ N^i(\bphi,\partial \bphi)(t,\vec{y}) \bphi^+_i(t,\vec{y})
= \int_{y^0=t}d\vec{y}\ N^i(\bphi,\partial \bphi)(y)\bphi^+_i(y),
\]
where we assume a summation over the repeated index $i$.
The expression (\ref{main-quantity}) can be written as $U_{t_1}^{t_2}(\Theta_{t_2}u)$, where
\[
 U_{t_1}^{t_2} = T\hbox{exp}\int_{t_1}^{t_2}dy^0\int_{\mathbb{R}^n}d\vec{y}\
 N^i(\bphi,\partial \bphi)(y)\bphi^+_i(y).
\]
We can then recover an expansion of this integral with terms analogous 
by using Wick's theorem with the commutation rules
\[
\left[\bphi^i(x), \bphi^j(y)\right] = 
 \left[\bphi^+_i(x), \bphi^+_j(y)\right] = 0,
\quad
\left[\bphi^+_i(x),\bphi^j(y)\right] = G^j_i(x,y) = G^j_{ix}(y),
\]
where $G^j_i:= (e^j,Ge_i)$ (here $(e^1,\cdots,e^N)$ is the dual basis of $(e_1,\cdots ,e_N)$).
In other words the $\bphi^+_i(x)$'s play the role of \emph{annihilation} operators and the 
$\bphi^i(x)$'s play the role of \emph{creation} operators.\\

As an example, 
we consider solutions $u$ of the scalar equation $\square_gu + u^3=0$ on a 4-dimensional space-time 
$(\mathcal{M},g)$ (see Section \ref{section-curved}) and we are given a smooth family of admissible Cauchy hypersurfaces
$(\sigma_s)_{s\in \R}$ which, for simplicity, we assume to be the level sets $\sigma_s = \tau^{-1}(s)$
of a temporal function $\tau\in \mathcal{C}^\infty_c(\mathcal{M})$.
We let $V(s,\varphi):= \Phi_{\sigma_s}(0,\lambda_s\varphi^3|_{\sigma_s})$ be the associated 
family of vector fields. We can express it more intuitively by setting
\[
 V_s = \int_{\sigma_s}d\mu_g(y)\lambda_s(y)\bphi(y)^3G_y,\quad\hbox{so that}\quad
 V_s(\varphi) = \int_{\sigma_s}d\mu_g(y)\lambda_s(y)\varphi(y)^3G_y.
\]
Then the corresponding first order operator reads
\[
 V_s\cdot = \int_{\sigma_s}d\mu_g(y)\lambda_s(y)\bphi(y)^3\bphi^+(y)
 = \int_{\sigma_s}d\overline{y}\bphi(y)^3\bphi^+(y),
\]
where we introduced the shorter notation $d\overline{y}:= d\mu_g(y)\lambda_s(y)$.
Let $f\in (\mathcal{E}_0^s)^*$ be linear, of the form $f = \int_\mathcal{M}d\hbox{vol}_g(x)\alpha(x)\bphi(x)$
(or equivalentely
$f(\varphi) = \int_\mathcal{M}d\hbox{vol}_g(x)\alpha(x)\varphi(x)$, $\forall \varphi\in \mathcal{E}_0^s$),
where $\alpha\in \mathcal{C}^\infty_c(\mathcal{M})$. Then
\[
 V_s\cdot f = \int_\mathcal{M}d\hbox{vol}_g(x)\alpha(x)
 \int_{\sigma_s}d\overline{y}  G_{y}(x) \bphi(y)^3
\]
and, writing $d\hbox{vol}_g(x) \simeq dx$ for short,
\[
 V_{s_2}\cdot(V_{s_1}\cdot f) = 3\int_\mathcal{M}\alpha(x)dx
 \int_{\sigma_{s_2}}d\overline{y_2} \int_{\sigma_{s_1}}d\overline{y_1}
 G_{y_1}(x)G_{y_2}(y_1)\bphi(y_2)^3\bphi(y_1)^2.
\]
We thus deduce the first terms in the expansion of $U_{t_1}^{t_2}f$ (relating the Cauchy data $\sigma_{t_1}$
and $\sigma_{t_2}$).
\[
 \begin{array}{ccl}
  U_{t_1}^{t_2}f & = & \displaystyle f + \int_{t_1}^{t_2}ds V_s\cdot f +
 \int_{t_1}^{t_2}ds_2\int_{t_1}^{s_2}ds_1V_{s_2}\cdot(V_{s_1}\cdot f) + \cdots\\
 & = & \displaystyle \int_\mathcal{M}\alpha(x)dx\,\bphi(x) + \int_{t_1}^{t_2}ds \int_\mathcal{M}\alpha(x)dx
 \int_{\sigma_s}d\overline{y}G_y(x)\bphi(y)^3 \\
 & & \displaystyle + 3\int_{t_1}^{t_2}ds_2\int_{t_1}^{s_2}ds_1 \int_\mathcal{M}\alpha(x)dx
\int_{\sigma_{s_2}}d\overline{y_2} \int_{\sigma_{s_1}}d\overline{y_1} G_{y_1}(x)G_{y_2}(y_1)\bphi(y_2)^3\bphi(y_1)^2
+ \cdots
 \end{array}
\]
Using $dsd\overline{y} = dsd\mu_g(y)\lambda_{s}(y) = d\hbox{vol}_g(y) \simeq dy$ and
setting $\int_{\sigma_{t_1}}^{\sigma_{t_2}}dy = \int_{t_1<\tau(y)<t_2}dy$,
\[
 \begin{array}{ccl}
 U_{t_1}^{t_2}f  & = & \displaystyle \int_\mathcal{M}\alpha(x)dx\,\bphi(x)
 + \int_\mathcal{M}\alpha(x)dx\int_{\sigma_{t_1}}^{\sigma_{t_2}}dy\,G_y(x)\bphi(y)^3 \\
 & &\displaystyle  + 3\int_\mathcal{M}\alpha(x)dx\,\int_{\sigma_{t_1}}^{\sigma_{t_2}}dy_2\int_{\sigma_{t_1}}^{\sigma_{\tau(y_2)}}dy_1
 G_{y_1}(x)G_{y_2}(y_1)\bphi(y_2)^3\bphi(y_1)^2 + \cdots.
 \end{array}
\]
Now apply Theorem \ref{theo-intro}: for any solution $u$ of $\square_gu+u^3=0$ and for $\|u\|$ and
$|t_2-t_1|$ sufficiently small, we have $f(\Theta_{t_1}u) = (U_{t_1}^{t_2}f)(\Theta_{t_2}u)$. Hence
assuming for simplicity that $u$ is continuous, $f = \bphi(x)$ for some $x\in \mathcal{M}$ and
$t_1>t_2$, we get 
\begin{equation}\label{formuleexemple}
 \begin{array}{c}
  \displaystyle \Theta_{t_1}u(x) = \Theta_{t_2}u(x)
 - \int_{\sigma_{t_2}}^{\sigma_{t_1}}dy\,G_y(x)(\Theta_{t_2}u(y))^3 \hfill \\
 \displaystyle \hfill +\, 3\int_{\sigma_{t_2}}^{\sigma_{t_1}}dy_2\int_{\sigma_{\tau(y_2)}}^{\sigma_{t_1}}dy_1
 G_{y_1}(x)G_{y_2}(y_1)(\Theta_{t_2}u(y_2))^3(\Theta_{t_2}u(y_1))^2 +\; \cdots 
 \end{array}
\end{equation}
Each term of the form $\Theta_tu(x)$ (marked by a line with a bold foot on the diagram below) reads:
\[
 \Theta_tu(x) = \int_{\sigma_t}d\mu_g(y)\left(\langle N,\nabla G_x\rangle_g(y)u(y) - G_x(y)\langle N,\nabla u\rangle_g(y)\right).
\]
Identity (\ref{formuleexemple}) (for e.g. $t_2<t_1\leq x^0$) is pictured by the following diagram representation
\begin{figure}[h]
\begin{center}
\input{diagram.pstex_t}
\end{center}
\end{figure}
where Feynman rules are used (see \cite{harrivel-helein}). Note that
if $x^0=t_1$, the l.h.s. of (\ref{formuleexemple}) is nothing but $u(x)$.

\section{A list of examples}
\textbf{Klein--Gordon equations}\\
The Main Theorem can be applied to all nonlinear Klein--Gordon equations of the type 
$\square u + m^2 + N(u) = 0$, where $u$ is a (real-valued) scalar field and $N$ is a
real analytic function (e.g. any polynomial or trigonometric function) for $s>n/2$.
However as already stressed in Remark \ref{remarksursegal1} this result extends straightforwardly
to the case $s=1\leq n/2$, if $N$ is a polynomial of degree
less than or equal to $n/n-2$.\\

\noindent
\textbf{Schr{\"o}dinger equations}\\
Our result can be applied only in the case where $n=1$, for any real analytic nonlinear function
$N$, i.e. to the equation $i\partial_0u + (\partial_1)^2u + N(u) = 0$ and for $s>1/2$, since
$H^s(\mathbb{R})$ is then an algebra and because of the continuous embedding $H^s(\mathbb{R})\hookrightarrow L^2(\mathbb{R})$.\\

\noindent
\textbf{Wave maps}\\
We consider for instance wave maps into the unit sphere $S^k\subset \mathbb{R}^{k+1}$ (but we
may replace $S^k$ by any Riemannian manifold which admits a real analytic isometric embedding in
some Euclidean space). We set $H^s(\mathbb{R}^n,S^k)
:= \{v\in H^s(\mathbb{R}^n,\mathbb{R}^{k+1}); v(x)\in S^k\hbox{ a.e.}\}$.
Wave maps are maps $u\in \mathcal{C}^0(\mathbb{R},H^s(\mathbb{R}^n,S^k))
\cap \mathcal{C}^1(\mathbb{R},H^{s-1}(\mathbb{R}^n,\mathbb{R}^{k+1}))$, which are weak solutions of
the system:
\[
 \square u + (|\partial_0u|^2 - |\vec{\partial}u|^2)u = 0.
\]
We note that the nonlinearity $N(u,\partial u) = (|\partial_0u|^2 - |\vec{\partial}u|^2)u$
is quadratic in $\partial u$ and hence does not satisfy (\ref{affine-hypothesis}). Thus our
result applies with $s>n/2+1$.\\

\noindent
\textbf{The Dirac--Maxwell system}\\
Set $n=4$, $g_{\mu\nu}:= \hbox{diag}(1,-1,-1,-1) = g^{\mu\nu}$ and consider $4\times 4$
Dirac matrices $\gamma^0,\gamma^1,\gamma^2,\gamma^3$ satisfying the Clifford algebra condition $\gamma^\mu\gamma^\nu +
\gamma^\nu\gamma^\mu = 2 g^{\mu\nu}$. We agree to sum over any repeated index.
The Dirac operator is $\dirac = \gamma^\mu\partial_\mu$,
acting on functions $\psi:\mathbb{R}^4\longrightarrow \mathbb{C}^4$. The Dirac--Maxwell system can be written
\begin{equation}\label{dirac-maxwell}
\left\{\begin{array}{ccc}
i\dirac \psi - m \psi & = & e\gamma^\mu\psi A_\mu\\
       \partial_{\nu} F^{\nu\mu} & = & e\overline{\psi}\gamma^\mu\psi,
       \end{array}
\right.
\end{equation}
where $A = A_\mu dx^\mu$ is a gauge connection for the electromagnetic field, $F_{\mu\nu}:= \partial_\mu A_\nu
- \partial_\nu A_\mu$ is the electromagnetic field and $F^{\mu\nu}:= g^{\mu\lambda}g^{\nu\sigma}F_{\lambda\sigma}$.
If we further assume the Lorentz gauge condition $\partial_\mu A^\mu = 0$, where $A^\mu:= g^{\mu\nu}A_\nu$, then
$\partial_\nu F^{\nu\mu} = \square A^\mu$, so that (\ref{dirac-maxwell}) can be written:
\begin{equation}\label{dirac-maxwell-var}
\left\{\begin{array}{ccc}
i\dirac \psi - m \psi - e\gamma^\mu\psi A_\mu  & = & 0\\
       \square A^\mu -e\overline{\psi}\gamma^\mu\psi & = & 0,
       \end{array}
\right.
\end{equation}
i.e. has the form (\ref{L+N-a}). Our result can hence be applied if
we assume that $A\in C^0(\mathbb{R},H^s(\mathbb{R}^3,(\mathbb{R}^4)^*)\cap C^1(\mathbb{R},H^{s-1}(\mathbb{R}^3,(\mathbb{R}^4)^*)$
and $\psi\in C^0(\mathbb{R},H^s(\mathbb{R}^3,\mathbb{C}^4)\cap C^1(\mathbb{R},H^{s-1}(\mathbb{R}^3,\mathbb{C}^4)$,
for $s>3/2$. Note that the Lorentz gauge can be achieved by starting from any arbitrary gauge connection
$\tilde{A}_\mu$ by setting $A_\mu = \tilde{A}_\mu - \partial_\mu \varphi$, where $\varphi(x) =
\int_0^{x^0}dy^0\int_{\mathbb{R}^3}d\vec{y}G(x-y)(\partial_\mu\tilde{A}^\mu)(y)$ (so that
$\square \varphi = \partial_\mu\tilde{A}^\mu$).\\

\noindent
\textbf{The pure Yang--Mills equation}\\
Given a finite dimensional semi-simple Lie algebra $\mathfrak{g}$,
the Yang--Mills equation for a connection $A:\mathbb{R}\times \mathbb{R}^n\longrightarrow \mathfrak{g}$
reads:
\begin{equation}\label{yang-mills}
 \partial_\nu\left(\partial^\nu A^\mu - \partial^\mu A^\nu + [A^\nu,A^\mu]\right)
+ \left[A_\nu,\partial^\nu A^\mu - \partial^\mu A^\nu + [A^\nu,A^\mu]\right] = 0,
\end{equation}
where we use the same convention on repeated indices as in the previous paragraph. Assuming again
that the Lorentz gauge condition $\partial_\mu A^\mu = 0$ is satisfied (which, in this nonlinear case, is
harder to achieve than for electromagnetism), then the higher order term is simply
$\partial_\nu\left(\partial^\nu A^\mu - \partial^\mu A^\nu\right) = \square A^\mu$. Then
(\ref{yang-mills}) has the form 
\begin{equation}\label{yang-mills-a}
 \square A^\mu + \underbrace{\left[A_\nu,[A^\nu,A^\mu]\right]}_{\hbox{cubic in }A}
+ \underbrace{ \partial_\nu\left([A^\nu,A^\mu]\right)+ \left[A_\nu,\partial^\nu A^\mu - \partial^\mu A^\nu \right]}_{\hbox{linear in }\partial A,
\hbox{ linear in }A}
= 0
\end{equation}
and hence satisfies Hypothesis (\ref{affine-hypothesis}). Thus Theorems
\ref{theo-zero} and \ref{theo-intro} can be applied to solutions
of (\ref{yang-mills-a}) if $s>n/2>s-r$, i.e., since $n=3$ and $r=1$, if
$3/2<s<5/2$.

\end{document}